\documentclass[smallextended,referee,envcountsect,]{svjour3}
\smartqed
\usepackage{graphicx}
\usepackage{amsmath,amsfonts,amssymb,amscd}
\usepackage{color}
\usepackage{geometry}
\usepackage{tikz}
\usepackage{graphics}
\usetikzlibrary{arrows}
\usepackage{tikz-er2}
\usetikzlibrary{positioning}
\usetikzlibrary{shadows}

\def\disp{\displaystyle}
\def\R{\mathbb R}
\def\ox{\bar{x}}
\def\ve{\varepsilon}
\def\B{\mathbb B}
\def\E{\widetilde E}
\def\dom{\mbox{\rm dom}\,}
\def\st{\stackrel}
\def\ov{\bar{v}}
\def\ra{\rangle}
\def\la{\langle}

\def\N{\mathbb{N}}

\def\gph{\mbox{\rm gph}\,}
\def\oR{\overline{\R}}
\def\kk{\kappa}

\def\gg{\gamma}
\def\Lm{\Lambda}

\journalname{}
\begin{document}

\title{Tilt Stability for Nonlinear Programs under Relaxed Constant Rank Constraint Qualification}

\author{Nguyen Huy  Chieu\and Nguyen Thi Quynh Trang \and   Nguyen Thi Hai Yen}

\institute{Nguyen Huy  Chieu  \at
             Department of Mathematics,  Vinh University
               Nghe  An, Vietnam, \\
              chieunh@vinhuni.edu.vn  \and
              Nguyen Thi Quynh Trang  \at
             Department of Mathematics,  Vinh University,
               Nghe  An, Vietnam\\
              quynhtrang@vinhuni.edu.vn
                                            \and
             Nguyen Thi Hai Yen, Corresponding author \at
             Faculty of Mathematics and Information Technology, The University of Da Nang - University of Science and Education,
               Da Nang, Vietnam\\
            nthyen$\_$kt@ued.udn.vn
    }

\date{Received: date / Accepted: date}

\maketitle

\begin{abstract}This paper investigates the tilt stability of local minimizers for nonlinear programs  under the relaxed constant rank constraint qualification in finite dimensions.  By employing a neighborhood primal-dual approach and extending calculus rules for subgradient graphical derivative,  we obtain some pointbased characterizations of tilt-stable local minimizers along with   an explicit formula for calculating the exact bound of tilt stability. These results extend the corresponding ones of H.~Gfrerer and B.~S.~Mordukhovich [SIAM J. Optim. 25 (2015), 2081-2119] by relaxing the constraint qualification and removing the linear independence condition of gradients of equality constraint functions. Examples are provided illustrating  our findings.
 \end{abstract}
\keywords{Tilt stability\and Pointbased second-order characterization\and Subgradient graphical derivative\and Nonlinear program \and Relaxed constant rank constraint qualification. }
\subclass{49J53\and  90C31\and  90C46}

\section{Introduction}

Tilt stability, which was introduced by Poliquin and Rockafellar \cite{PR98}, is widely recognized as a crucial property of  local minimizers.  It guarantees the local minimizer evolves uniquely in  a Lipschitz continuous trajectory  under small linear perturbations. This behavior is highly desirable for local minimizers, offering benefits for both theoretical analysis and numerical computations.   Tilt stability is closely related to some other significant concepts, including strong variational convexity, uniform second-order growth condition, and strong metric regularity of the subdifferential \cite{BS,DL13,KMP23,MN1}. Due to its importance and connections, tilt stability has attracted considerable attention from many researchers over the years \cite{BGM19,CHN18,CHT20,DMN14,EW12,G25,KMP23,KMPV25,LZ13,M24,MN13,MN1,MO13,MTW25,MS21}.

To determine if a stationary point of an optimization problem is a tilt-stable local minimizer, among other things, researchers often employ second-order generalized differentiations and exploit the specific structural properties of the problem. This method has yielded   various second-order characterizations of tilt stability, which can be  classified  into neighborhood and pointbased types~\cite{M24}. While neighborhood characterizations involve points near the reference point, pointbased characterizations are formulated entirely at the local minimizer. The latter is generally more desirable for practical applications.

The first neighborhood second-order characterization of tilt stability for unconstrained optimization problems with continuously prox-regular extended-real-valued objectives was established by Mordukhovich and Nghia \cite[Theorem 3.5]{MN1} via the combined second-order subdifferential. Subsequently,  applying this to classical nonlinear programs (NLPs) and utilizing calculus rules for the combined second-order subdifferential, they proved in \cite[Theorem 4.3]{MN1} that the uniform second-order sufficient condition (USOSC) serves as a neighborhood characterization of tilt-stable local minimizers for NLPs under both the Mangasarian-Fromovitz constraint qualification (MFCQ) and the constant rank constraint qualification (CRCQ).
Later, Chieu~et~al.~\cite{CHN18} provided an alternative neighborhood characterization of tilt stability for unconstrained optimization problems with continuously prox-regular extended-real-valued objectives. This characterization relies on the subgradient graphical derivative.  They then used it to show in \cite[Theorem 4.5]{CHN18} that the relaxed uniform second-order sufficient condition (RUSOSC), which reduces  to USOSC under CRCQ,  is  a neighborhood characterization of tilt-stable local minimizers for NLPs satisfying the metric subregularity constraint qualification (MSCQ).

  The   pointbased second-order characterization of tilt stability for uncontrained optimization problems with continuously prox-regular extended-real-valued objectives was  discovered  by Poliquin and Rockafellar \cite[Theorem 1.3]{PR98}. It is formulated as  the positive semi-definiteness of   the limiting second-order subdifferential of the objective function at a stationary point relative to the origin. The result was then applied  to characterizing tilt stability for  the NLPs with linear constraints \cite[Theorem 4.5]{PR98} and  the NLPs satisfying  the  linear independent constraint qualification (LICQ) \cite[Theorem~5.2]{MR}.
    This pointbased dual-dual approach to characterizing tilt stability for NLPs   requires the exact computation of the limiting second-order subdifferential, which  is often a challenging task, especially  under weak constraint qualifications. Using Mordukhovich and Nghia's neighborhood   characterization of tilt stability \cite[Theorem~3.5]{MN1},   Gfrerer and Mordukhovich established pointbased second-order characterizations of tilt stability together with explicit formulas calculating  the exact bound of tilt stability for the NLPs under either nondegeneracy or 2-regularity \cite[Theorem 7.6]{GM15} and  under CRCQ   \cite[Theorem 7.7]{GM15}.
    This neighborhood dual-dual approach hinges on the boundedness of a specific sequence of multipliers, which is ensured by the bounded extreme point property \cite[Definition~3.3]{GM15}.  Besides, Gfrerer and Mordukhovich also demonstrated  \cite[Example~8.5]{GM15} that a pointbased second-order characterization of tilt stability for NLPs under MFCQ alone  is  impossible. As a result, there is no hope to have a pointbased second-order characterization of tilt stability for NLPs  if we rely solely on MFCQ or a constraint qualification that is weaker than MFCQ.

The relaxed constant rank constraint qualification (RCRCQ) for NLPs, proposed by Michenko and Stakhovski \cite{MS11}, is among significant constraint qualifications  that are weaker than CRCQ.
While RCRCQ requires checking fewer constraints for satisfaction compared to CRCQ, it shares many common features with CRCQ \cite{L12,ML11b,RS24}. In particular,  RCRCQ is stronger than MSCQ and independent of MFCQ. Like CRCQ,   RCRCQ is also a second-order constraint qualification  for  the strong second-order necessary optimality condition \cite{ML16,RS24}.
Despite the crucial role of tilt stability in nonlinear programming and the appealing characteristics of RCRCQ, pointbased second-order characterizations of tilt stability for the NLPs under the RCRCQ remain uninvestigated. Our main goal in this paper is  to fill this gap by addressing the natural question:  Can we derive pointbased second-order characterizations of tilt-stable local minimizers for NLPs under RCRCQ? We also want to determine  if the linear independence condition for gradients of equality constraint functions can be removed from the assumptions of    \cite[Theorem~7.7]{GM15}.

To achieve our goal, we adopt an approach distinct from Gfrerer and Mordukhovich's neighborhood dual-dual method in \cite[Theorem 7.7]{GM15}.  We instead utilize the neighborhood second-order characterization of tilt stability via the subgradient graphical derivative, as developed by Chieu et al. \cite{CHN18}.
This neighborhood primal-dual approach together with extending calculus rules for subgradient graphical derivative
allows us to obtain pointbased characterizations of tilt-stable local minimizers along with  an explicit formula for calculating the exact bound of tilt stability for NLPs under RCRCQ in the absence of  the linear independence of   the  gradients of equality constraint functions. Therefore, the answers to the aforementioned questions are affirmative.

The rest of the paper is organized as follows. Section 2 reviews some  concepts and their properties from variational analysis and nonlinear programming, which are needed for our analysis. We  specifically recall the concept of tilt-stable minimizers and its neighborhood second-order characterization via the subgradient graphical derivative \cite[Theorem 3.3]{CHN18}. While the second-order subdifferentials used in \cite{GM15,MN1,PR98} are considered as dual-dual second-order generalized differential constructions, the subgradient graphical derivative is a primal-dual one. This section also covers the classical first-order necessary optimality condition and significant constraint qualifications for the NLPs, along with their relationships.
Section 3 analyzes the graphical derivative of normal cone mappings associated with the feasible sets of NLPs under RCRCQ.
 We show that the formula for this derivative, previously established under the weaker MSCQ \cite[Theorem~3.5]{CH17}, can be simplified.
  Section~4 focuses on characterizing the tilt stability and calculating the exact bound of tilt stability for the NLPs under RCRCQ.
 The key contribution here is Theorem~4.1, which extends \cite[Theorem 7.7]{GM15} by achieving the same results under weaker, more verifiable assumptions. Specifically, we relax the CRCQ to RCRCQ and eliminate the requirement for linear independence of the gradients of equality constraint functions. This section also includes some
  examples illustrating  our findings. Section 5 concludes the paper with a summary of our findings and some perspectives for future work.

\section{Preliminaries}

This section reviews essential concepts and their properties from variational analysis and nonlinear programming \cite{AES10,AHSS11,AHSS11b,AMS05,KMO14,M1,M24,RW98}, which are fundamental to our analysis.

Unless otherwise specified, this paper uses the following standard notations:
 $\R^n$   is   a Euclidean space with the inner product $\langle\cdot,\cdot\rangle$ and norm $\|\cdot\|$, and we  identify this space with its dual;    the closed ball centered at $\ox$ with radius $\ve>0$ is denoted by $\B_\ve(\ox)$, and   $\B_{\R^n}$ is the closed unit ball in $\R^n$;
    for any finite set $J$, $|J|$  represents the number of elements in $J$;  the rank of a family of vectors $\{a_j\, |\, j\in J\}$ is  ${\rm rank}\{a_j\ |\ j\in J\}$;
 $\overline{\R}:=\R\cup\{\infty\}$, and $\N:=\{0,1,2,...\}$.

 Let    $f: \R^n\rightarrow \overline{\R}$  be   a lower semicontinuous proper function.
Recall  \cite{M1,M24,RW98} that the {\it limiting  subdifferential}   of $f$ at $\bar x\in \dom f:=
   \{x\in \R^n|\; f(x)<\infty\}$ is the set $\partial f(\bar x)$ in the dual space $\R^n$ defined by
		\begin{equation*}\label{ls} \partial f(\bar x):=\left\{ x^*\in\R^n\, |\, \exists x_k\to \bar x, x_k^*\to x^* \text{ with } f(x_k)\to f(\bar x),\,  x_k^*\in \widehat{\partial}f(x_k), \, k\in \N \right\}, \end{equation*}
where
		 \begin{equation*}\label{rs} \widehat{\partial}f(x_k):=\left\{ x^*\in\R^n\, \Big|\, \liminf\limits_{x\rightarrow x_k}\frac{f(x)-f(x_k)-\la x^*,x-x_k\ra}{\|x- x_k\|}\geq 0\right\}\end{equation*}
is called the {\it regular subdifferential} of $f$ at $x_k$.   If $|f(\bar x)|=\infty$,    $\partial f(\bar x)=\widehat{\partial} f(\bar x):=\emptyset$ by convention.

	For a convex function f, the limiting and regular subdifferentials coincide with the standard subdifferential from convex analysis:
$$\partial f(\bar x)=\widehat{\partial} f(\bar x)=\begin{cases}\left\{x^*\in \R^n\, |\, f(x)-f(\bar x)\geq \la x^*, x-\bar x\ra\quad  \forall x\in \R^n\right\}\quad &\mbox{if}\ \bar x\in \dom f,\\
\quad\quad\quad\quad\quad\quad\quad\quad\emptyset &\mbox{otherwise.}\end{cases}$$
 Furthermore, if  $f$ is  continuously differentiable at $\bar x$, we have
 $$\partial f(\bar x)=\widehat{\partial} f(\bar x)=\{\nabla f( \bar x)\},$$ where $\nabla f( \bar x):=\left(\frac{\partial f}{\partial x_1}(\bar x),\frac{\partial f}{\partial x_2}(\bar x),....,\frac{\partial f}{\partial x_n}(\bar x)\right).$

For a set $\mathrm{\Omega}\subset\R^n,$ the (Bouligand-Severi) {\it tangent cone} to  $\mathrm{\Omega}$ at $\bar x\in \mathrm{\Omega}$  is defined by
    \begin{equation*}\label{tcone}
  	T_{\mathrm{\Omega}}(\bar x):=\big\{u\in\R^n | \, \exists\  t_k \downarrow 0, \  u_k\rightarrow u\ \mbox{ with }\  \bar x+t_ku_k\in\mathrm{\Omega} \ \ \forall k\in \N\big\},
\end{equation*}
	and the  {\it limiting normal cone} to $\mathrm{\Omega}$ at $\bar x\in \mathrm{\Omega}$ is the set $N_{\mathrm{\Omega}}(\bar x)$ in the dual space $\R^n$  given by
$$N_{\mathrm{\Omega}}(\bar x)=\left\{x^*\in \R^n\, | \,  \exists x_k\rightarrow \bar x,\, x_k^*\rightarrow x^*  \  \mbox{with} \  x_k^*\in \widehat{N}_\mathrm{\Omega}(x_k)\;\  \forall   k\in \N\right\},$$
where
$$\widehat{N}_{\mathrm{\Omega}}(x_k):=\left\{x^*\in\R^n\, \Big|\, \limsup\limits_{x\st{\mathrm{\Omega}}\rightarrow x_k}\frac{\langle x^*, x-x_k\rangle}{\|x-x_k\|}\leq 0\right\}$$
is the  {\it regular normal cone} to $\mathrm{\Omega}$ at $x_k \in \mathrm{\Omega}$. If $\bar x\in \R^n\backslash\mathrm{\Omega}$, $N_{\mathrm{\Omega}}(\bar x)=\widehat{N}_{\mathrm{\Omega}}(\bar x):=\emptyset$ by convention.

It is known \cite[Theorem~1.10]{M1}  that the regular normal cone is the (negative) dual  of the tangent cone:
  \begin{equation*}\label{eqdualTN}\widehat{N}_{\mathrm{\Omega}}(\bar z)=\big(T_{\mathrm{\Omega}}(\bar z)\big)^*,
  \end{equation*}
  where $C^*:=\{x^*\in \R^n\ |\ \la x^*,u\ra\le 0\quad \forall u\in C\}$ for any $C\subset \R^n$.
 If $\Omega$ is convex, the above tangent and normal cones  reduce to the tangent  and normal cones in the sense of convex analysis, respectively.
Observe that
 $$ \widehat{\partial}\delta_{\mathrm{\Omega}}(x) =  \widehat{N}_{\mathrm{\Omega}}(x) \text{ and } \partial \delta_{\Omega}(x) =  N_{\mathrm{\Omega}}(x) \text{ for all } x\in \mathrm{\Omega},$$
  where $\delta_{\mathrm{\Omega}}(x)$ stands for the indicator function of the set $\mathrm{\Omega}$, which is equal
 to 0 if $x\in\mathrm{\Omega}$ and to $\infty$ otherwise.

  Let $\mathrm{\Phi}: \R^n\rightrightarrows\R^m$ be a set-valued mapping with its graph  ${\rm gph}\mathrm{\Phi}:=\big\{(x,y)\,  |\,   y\in \mathrm{\Phi}(x)\big\}$ and  its domain ${\rm dom} \mathrm{\Phi}:=\big\{x\,  |\,  \mathrm{\Phi}(x)\not=\emptyset\big\}.$
The {\it graphical derivative}  of $\mathrm{\Phi}$ at a given point  $(\bar x,\bar y)\in{\rm gph}\mathrm{\Phi}$ is the set-valued mapping $D\mathrm{\Phi}(\bar x,\bar y):\R^n\rightrightarrows \R^m$ defined by
  \begin{equation*}
 D\mathrm{\Phi}(\bar x,\bar y)(u):=\big\{ v\in\R^m\, |\, (u,v)\in T_{{\rm gph}\mathrm{\Phi}}(\bar x,\bar y)\big\}\ \,  \mbox{for all}  \ u\in\R^n.
 \end{equation*}
   In the case  $\mathrm{\Phi}(\bar x)=\{\bar y\},$  one writes $D\mathrm{\Phi}(\bar x)$  for $D\mathrm{\Phi}(\bar x,\bar y)$. If $\mathrm{\Phi}$ is a
   single-valued mapping being Fr\'echet differentiable at $\bar x$, then
   $$ D\mathrm{\Phi}(\bar x)(u)= \{\nabla \mathrm{\Phi}(\bar x)u\} \text{ for all } u\in \R^n,$$
where $\nabla \mathrm{\Phi}(\bar x)$ is the Jacobian matrix of $f$ at $\bar x$.

Let $f: \R^n\rightarrow \overline{\R}$ be   a lower semicontinuous function.  Recall  \cite[Chapter 13]{RW98} that $f$ is  said to be {\it prox-regular} at $\ox\in \dom f$ for $\ov^*\in \partial f(\ox)$ if there exist real numbers $r,\ve>0$ such that for all $x,u\in \B_\ve(\ox)$ with $|f(u)-f(\ox)|\leqslant\ve$ we have
\begin{equation*}\label{prox}
	f(x)\ge f(u)+\la v^*, x-u\ra-\frac{r}{2}\|x-u\|^2\quad \mbox{for all}\quad v^*\in \partial f(u)\cap \B_\ve(\ov^*);
\end{equation*}
and $f$  is called {\em subdifferentially continuous} at $\ox$ for $\ov^*\in \partial f(\ox)$ if $f(x_k)\to f(\ox)$  whenever $x_k\to \ox$ and  $v_k^*\to \ov^*$  with $v_k^*\in \partial f(x_k)$, $k\in \N$.
	
	If $f$ is subdifferentially continuous at
	at $\ox$ for $\ov$, the inequality ``$|f(u)-f(\ox)|\leqslant \ve$'' in the definition of prox-regularity above
	could be omitted. When $f$ is both prox-regular and subdifferentially continuous at
	at $\ox$ for $\ov$, we say that $f$ is {\it continuously prox-regular}  at $\ox$ for $\ov$.

 \begin{definition}{ \rm (\cite{M24,PR98}) \label{tilt}  Let $f\colon\R^n\to\oR$, and let $\ox\in\dom f$.
			
	$(i)$ The point $\bar x$ is said to be  a  {\it tilt-stable local minimizer} of $f$  if there is a number $\gg>0$ such that the mapping
	\begin{eqnarray}\label{2.15}
		M_\gg(v):={\rm argmin}\big\{f(x)-\la v,x\ra\big|\;x\in \mathbb{B}_\gg(\ox)\big\}, \quad v\in\R^n
	\end{eqnarray}
	is single-valued and Lipschitz continuous on some neighborhood of $\bar v=0\in \R^n$ with $M_\gg(0)=\{\ox\}$.
	
	$(ii)$ Given $\kk>0$, the point $\bar x$ is called a  {\it tilt-stable local minimizer} of $f$ with {\sc modulus} $\kk$, if there is a number $\gg>0$ such that $M_\gg(0)=\{\ox\}$ and the mapping $M_\gg$  in \eqref{2.15}
	is single-valued and Lipschitz continuous with modulus $\kk$ around the origin $0\in \R^n$.
	
	$(iii)$ The {\it exact bound of tilt stability} of function $f$ at $\ox$ is defined by
	\[
	{\rm tilt}\, (f,\ox):=\inf_{\gg>0} {\rm lip} M_\gg(0) .
	\]
	via the exact Lipschitzian bound of the mapping $M_\gg$  from \eqref{2.15} around the origin.
}
\end{definition}

The following theorem provides  the  neighborhood second-order characterization of tilt-stable minimizers for extended real-valued functions, along with a precise formula for calculating the exact bound of tilt stability.
\begin{theorem}\label{theo22} {\rm (\cite[Theorem 3.3]{CHN18})}
	Let $f:\R^n\to\oR$ be a l.s.c.\ proper function with $\ox\in \dom f$ and  $0\in\partial f(\ox)$. Assume that $f$ is continuously prox-regular  at $\ox$ for $\ov=0$. Then the following assertions are equivalent:
	
	{\it (i)} The point $\ox$ is a tilt-stable local minimizer of $f$ with modulus $\kk>0$.
	
	{\it (ii)} There is a constant $\eta>0$ such that for all $w\in\R^n$ we have
	\begin{equation*}\label{3.2}
		\la z,w\ra\ge\frac{1}{\kk}\|w\|^2\;\mbox{ whenever }\; z\in \partial_D^2f(x,v)(w)\;\mbox{ with }\;(x,v)\in\gph\partial f\cap\B_\eta(\ox,0),
	\end{equation*}
where $\partial_D^2f(x,v):=D\partial f(x,v)$ is the subgradient graphical derivative of $f$ at $x$ relative to $v$.
	Furthermore, the exact bound of tilt stability of $f$ is calculated by the formula
	\begin{eqnarray*}\label{3.2b}
		{\rm tilt}\, (f,\ox)=\inf_{\eta>0}\sup\Big\{\frac{\|w\|^2}{\la z,w\ra}\Big|\; z\in {\partial}_D^2f (x, v)(w),\;(x, v)\in\gph\partial f\cap\B_\eta(\ox,0)\Big\}
	\end{eqnarray*}
	with the convention that $0/0:=0$.
\end{theorem}

In this paper, we consider the following classical nonlinear program:
 \begin{equation}\label{tiltMP1} \begin{array}{rl} &\mbox{minimize} \quad \varphi(x)\quad \mbox{subject to} \quad q_i(x)=0,\;i\in E \text{ and } q_i(x)\leq 0,\ i\in I,
 	\end{array}
 \end{equation}
 where the objective function $\varphi:\R^n\rightarrow\R$ and constraint functions $q_i: \R^n\to \R,$ $i\in E\cup I,$  are  twice continuously differentiable. Here  $E := \{i\in \N\ |\ 1\leq i\leq {\ell}_1\}$ and $I := \{i\in \N\ |\ {\ell}_1 + 1\leq i\leq {\ell}_1 + {\ell}_2\}$, with $\ell_1,\ell_2\in \N$ and $E\cup I\not=\emptyset$.

The {\it feasible set} $\Gamma$ of \eqref{tiltMP1} is defined as the solution set of  the following system of  equality and inequality constraints:
 \begin{equation}\label{CSy}
\begin{cases}
	q_i(x)=0& \text{ for } i\in E,\\
	q_i(x)\leq 0& \text{ for } i\in I.\\
\end{cases}
\end{equation}
In other words, the feasible set $\Gamma$ is the one given by
$$\Gamma:=\big\{x\in \R^n\ |\ q_i(x)=0\  \text{ for } i\in E,\ \mbox{and}\
	q_i(x)\leq 0\  \text{ for } i\in I\big\}.$$
Each point $x\in \Gamma$ is said to be a {\it feasible point} (or {\it feasible solution}) of problem \eqref{tiltMP1}.
If $\ell_1=0$ then $E=\emptyset$ and the system \eqref{CSy} has no equality constraint while if $\ell_2=0$ then $I=\emptyset$ and  the system \eqref{CSy} has no inequality constraint.

Let $\;\Theta:=\{0_{\R^{{\ell}_1}}\}\times \mathbb{R}_{-}^{{\ell}_2}$ and  let $q: \R^n\rightarrow\R^{\ell}$ be the twice continuously differentiable mapping given by $q(x)=\big(q_1(x), \ldots , q_{\ell}(x)\big)$  with ${\ell} := {\ell}_1 + {\ell}_2.$
Then
 \begin{equation}\label{CS}\Gamma=\{x\in \R^n\, |\, q(x)\in\Theta\}.\end{equation}

The {\it linearized tangent cone} to $\Gamma$ at $x\in \Gamma$  (with respect to the representation \eqref{CS})
is the set  $T^{\rm lin}_{\Gamma}(x)$ defined by
 \begin{equation*}\label{ltc}
 	T^{\rm lin}_{\Gamma}(x):=\big\{u\in \R^n\;|\; \la \nabla q_i(x),u\ra =0\ \mbox{for} \; i\in E \text{ and } \la \nabla q_i(x),u\ra \leq 0\; \mbox{for}\  i\in I(x)\big\},
\end{equation*}	
where $I(x):=\{i\in I\;|\; q_i(x)=0\}$ is the {\it  index set of active inequality constraints}.

The function $\mathcal{L}: \R^n\times\R^{\ell}\to \R$ defined  by  $\mathcal{L}(x,\lambda)=\varphi(x)+\big\la\lambda,q\big\ra(x)$ is called the {\it Lagrange function} associated with problem \eqref{tiltMP1}, where
 $\big\la\lambda,q\big\ra(x):=\big\la\lambda,q(x)\big\ra$ for $x\in \R^n$ and $\lambda\in \R^\ell$.
The Jacobian matrix and Hessian matrix of the function $\mathcal{L}(\cdot,\lambda)$  at $\bar x$ are  denoted by  $\nabla_x\mathcal{L}(\ox, \lambda)$ and $\nabla_x^2\mathcal{L}(\ox, \lambda),$ respectively.

  Given  a point  $\bar x\in \Gamma$,  the {\it set of  Lagrange multipliers} of  problem \eqref{tiltMP1} at $\bar x$  is the one defined as
 \begin{equation*}
 \label{som1}
 	\Lambda(\bar x):=\left\{\lambda\in \R^{\ell_1}\times \R_+^{\ell_2}\ |\  \nabla_x\mathcal{L}(\ox, \lambda)=0\ \mbox{and}\    \lambda_iq_i(\bar x)=0\quad \mbox{for} \ i\in I\right\},
 \end{equation*}
 and   $\bar x$ is called  a {\it  stationary point} of problem \eqref{tiltMP1} if  $\Lambda(\bar x)\not=\emptyset.$

  By a constraint qualification, we mean a condition imposed on the feasible set and its description that guarantees the existence of Lagrange multipliers at a local minimizer.  It is worth emphasizing that the validity of a constraint qualification depends not only on the feasible set but also on its chosen representation.

Here are some well-known constraint qualifications in nonlinear programming \cite{AHSS11,GM15}
\begin{definition}{\rm Let $\ox$ be a feasible point of  problem \eqref{tiltMP1}. One  says that:
		
	$\bullet$ The {\it linear independence constraint qualification} (LICQ) holds at $\bar{x}$ if  $\{\nabla q_i(\bar x)|\; i\in E\cup {I}(\bar x)\}$ is  linearly independent.
	
	$\bullet$  The {\it  Mangasarian-Fromovitz constraint qualification} (MFCQ)  holds at $\bar{x}$ if  $\{\nabla q_i(\bar x)|\; i\in E\}$ is linearly independent and there exists $v\in \R^n$ such that $\la \nabla q_i(\bar x),v\ra=0$ for $i\in E$ and $\la \nabla q_i(\bar x),v\ra<0$ for $i\in I(\bar x)$.
			
	$\bullet$ The  {\it  constant rank constraint qualification} (CRCQ) holds at $\bar{x}$ if there exists a neighborhood
			$U$ of $\bar x$ such that for every $K\subseteq E$ and every $J\subseteq I(\bar x)$, the family of
			gradients $\{\nabla q_i(x)|\; i\in K\cup J\}$ has the same rank for every $x\in U$.
			
	$\bullet$  The {\it relaxed constant rank constraint qualification} (RCRCQ) holds at $\bar{x}$ if there exists a neighborhood
			$U$  of $\bar x$ such that for every $J\subseteq I(\bar x)$,  the family of gradients $\{\nabla q_i(x)|\; i\in E\cup J\}$ has the same rank for  every $x\in U$.	

$\bullet$ The {\it metric subregularity constraint qualification} (MSCQ) holds at $\bar x$ if  there exist a neighborhood $U$ of $\ox$ and a constant $\kk>0$ such that
		\begin{equation*}\label{subre}
			d(x;\Gamma)\le \kk d\big(q(x);\Theta\big)\quad \mbox{for all}\  x\in U.
		\end{equation*}

	$\bullet$ The {\it Abadie constraint qualification} (ACQ) holds at $\bar x$ if $T_\Gamma(\bar x) = T_\Gamma^{\rm lin}(\bar x).$
}
\end{definition}

The diagram below shows  some relationships among  the aforementioned constraint qualifications.  For a more detailed understanding of how constraint qualifications in nonlinear programming relate to each other and to numerical optimization issues, the readers can refer to \cite{AES10,AHSS11,AHSS11b,AMS05} and the references therein.
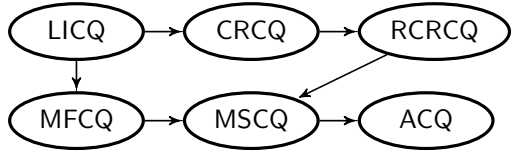
\begin{center}
\begin{tikzpicture}[
  font=\sffamily,
  every matrix/.style={ampersand replacement=\&,column sep=1cm,row sep=1cm},
  source/.style={draw,thick,rounded corners,top color=white, bottom color=green!20,
draw=green!50!black!100, drop shadow,inner sep=.2cm},
  process/.style={draw,thick,circle,fill=blue!20},
  sink/.style={source,top color=white, bottom color=red!20,
draw=red!50!black!100, drop shadow},
  datastore/.style={draw,very thick,shape=datastore,inner sep=.2cm},
  dots/.style={gray,scale=2},
  to/.style={->,>=stealth',shorten >=1pt,semithick,font=\sffamily\footnotesize},
  every node/.style={align=center}
   every attribute/.style  = {top color=white, bottom color=yellow!20,
draw=yellow, node distance=7em, drop shadow}]
     \node[attribute] (LI) {LICQ};
  \node[attribute] (MF) [below =0.4cm of LI] {MFCQ} edge (LI);
  \draw[to] (LI) -- (MF);
  \node[attribute](CR)[right=0.5 cm of LI] {CRCQ} edge (LI);
  \draw[to] (LI) -- (CR);
  \node[attribute] (RCR) [right=0.5 cm of CR] {{RCRCQ}} edge (CR);
  \draw[to] (CR) -- (RCR);
  \node[attribute] (MS) [right =0.5cm of MF] {MSCQ} edge (MF);
  \draw[to] (RCR) -- (MS);
   \draw[to] (MF) -- (MS);
  \node[attribute] (ACQ) [right =0.5cm of MS] {ACQ} edge (MS);
  \draw[to] (MS) -- (ACQ);
\end{tikzpicture}
\end{center}

\section{Graphical Derivative of  Normal Cone Mapping  under RCRCQ}
 Theorem 3.5 in \cite{CH17} gives us  a formula for computing graphical derivatives of  normal cone mappings associated with feasible sets of nonlinear programs  under MSCQ.   Our main aim in this section is to show that this formula can be simplified when  MSCQ is replaced by RCRCQ.

Let  $x\in \Gamma$ and $x^*, v\in\R^n$, where $\Gamma$ is defined by  \eqref{CS}. Recall \cite{CH17,GM15} that the {\it set of multipliers} associated with $(x, x^*)$ is given by
 \begin{equation}\label{som}
 \begin{array}{rl}	
 \Lambda(x, x^*)&:=\{\lambda\in N_{\Theta}\big(q(x)\big)\;|\; \nabla q(x)^*\lambda=x^*\}\\
 &=\{\lambda\in \R^{\ell}\ |\ \sum\limits_{i\in E\cup I}\lambda_i\nabla q_i(x)=x^*, \ \lambda_i\geq 0\ \mbox{and}\  \lambda_iq_i(x)=0 \  \mbox{as} \ i\in I\},
 \end{array}
 \end{equation}
 and the {\it critical cone} to $\Gamma$ at $(x,x^*)\in {\rm gph}{\partial}\delta_{\Gamma}$ is   defined  by
  $
	K(x,x^*):=T_\Gamma(x)\cap \{ x^*\}^\bot.
$
Put
  \begin{equation*}\label{msidp}
	\Lambda(x, x^*; v):={\rm argmax}\{ \la v, \nabla^2\la \lambda, q\ra (x)v\ra \big|\; \lambda\in\Lm(x, x^*)\}
\end{equation*}
and \begin{equation}\label{icong}I^+(x,x^*):=\bigcup_{\lambda\in\Lm(x,x^*)}I^+(\lambda),\end{equation}
where  \begin{equation}\label{sci} I^+(\lambda):=\{i\in I\;| \;\lambda_i>0\}\quad\mbox{for}\ \lambda\in\R^\ell.\end{equation}
 It follows from \eqref{som}   that   $q_i(x)=0$ whenever $i\in I^+(\lambda)$ and $\lambda\in \Lm(x,x^*)$.

The following result is taken from \cite{GM15}, which is needed for our subsequent analysis.
\begin{lemma}\label{cc4}{\rm (\cite[Proposition 4.3]{GM15})}
		Let MSCQ hold at $\bar x\in \Gamma$, and let $(x,x^*)\in {\rm gph}{\partial}\delta_{\Gamma}$ be any pair such that $x\in \Gamma$ is sufficiently close to $\ox$. Then the following assertions are satisfied:
		
		$(i)$ For every $\lambda\in \Lambda(x, x^*)$ we have
		\[	K(x,x^*)=\left\lbrace v\in\R^n\Big| \la \nabla q_i(x), v\ra\begin{cases}
			=0 & \text{ if } i\in E\cup I^+(\lambda) \\
			\leq 0 & \text{ if } i\in {I}(x)\setminus I^+(\lambda)
		\end{cases}\right\rbrace. \]

		$(ii)$ There exist  $\tilde{\lambda}\in \Lambda(x, x^*)$  and  $v\in K(x,x^*)$ such that   $I^+(\tilde{\lambda})=I^+(x,x^*)$ and
		\[ \la \nabla q_i(x), v\ra\begin{cases}
			=0 & \text{ if }  i\in E\cup I^+(\tilde{\lambda}), \\
			< 0 & \text{ if }  i\in {I}(x)\setminus I^+(\tilde{\lambda}).
		\end{cases}\]
	
		$(iii)$ For every   $v\in K(x,x^*)$ we have $\Lambda(x, x^*; v)\not=\emptyset$.
\end{lemma}

The following result was  justified by Gfrerer and Outrata  \cite[Lemma 4]{GO16} for the case where $E=\emptyset.$  Their proof can be applied to the case where $E\not=\emptyset$ with some minor modifications.
\begin{lemma}\label{GO16lm4}{\rm (\cite[Lemma 4]{GO16})}
Let $MSCQ$ hold at $x\in \Gamma$ and  $x^*\in N_\Gamma(x).$ Assume that
$$\Lm(x,x^*;v_1)=\Lm(x,x^*,v_2)\quad \mbox{whenever}\   v_1,v_2\in K(x,x^*)\backslash\{0\}.$$ Then
$$\nabla^2\la \lambda^1-\lambda^2,q\ra(x)v\in {\rm span}\big\{\nabla q_i(x)^*\ |\ i\in E\cup I^+(x,x^*)\big\}, $$
for every  $v\in K(x,x^*)$  and $\lambda^1,\lambda^2\in {\Lm}(x,x^*;v).$
\end{lemma}
The main result of this section reads as follows.
\begin{theorem} \label{prop53}
	Suppose that  RCRCQ is valid  at $\bar x\in \Gamma$. Then there exist an open neighborhood $U$ of $\bar x$ and a real number $\gamma>0$  such that MSCQ is fulfilled at every point $x\in \Gamma\cap U$ with modulus $\gamma>0,$ and
the following assertions hold:
\begin{itemize}
\item[$(i)$] For every $x\in \Gamma\cap U$
and $x^*, v\in\R^n$ satisfying
\begin{equation}\label{msid}
	 \la \nabla q_i(x), v\ra =0 \text{ whenever } i\in E\cup I^+(x,x^*),
\end{equation}
we have  $\Lambda(x, x^*; v)=\Lambda(x, x^*)$.
\item [$(ii)$]  For $x\in U\cap \Gamma,$   $x^*\in\R^n$ and  $\lambda\in \Lm(x,x^*)$ we have
  \begin{equation}\label{tc}
T_{{\rm gph}N_{\Gamma}}(x, x^*)=\big\{(w, z)\in\R^{2n}\big|\; z\in \nabla^2\la \lambda, q\ra (x)w+ N_{K(x,x^*)}(w)\big\},
\end{equation}
or equivalently,
$$
DN_\Gamma(x,x^*)(w)=\nabla^2\la \lambda, q\ra (x)w+ N_{K(x,x^*)}(w)\quad \mbox{for all}\ w\in \R^n.
$$
  \item [$(iii)$]  For $x\in U\cap \Gamma$ and  $x^*\in \R^n$ with $\Lambda(x, x^*)\not=\emptyset$ there exists $\lambda\in \Lambda(x, x^*)\cap\gamma\|x^*\|\B_{\R^\ell}$ such that
        $$z\in \nabla^2\la \lambda, q\ra (x)w+ N_{K(x,x^*)}(w)\quad\mbox{whenever}\ w\in \R^n\ \mbox{and}\ z\in DN_\Gamma(x,x^*)(w).$$
\end{itemize}
\end{theorem}
\noindent{\it Proof.} Since RCRCQ is fulfilled  at $\bar x$, there exists an open neighborhood $U$ of $\bar x$ such that $I(u)\subset I(\bar x)$ for every $u\in \Gamma\cap U$,  and
for each ${J}\subset I(\bar x)$ the system
$\{\nabla q_i(u)|\; i\in E\cup {J}\}$ has the same rank for every $u\in U$.  Furthermore, shrinking $U$ if necessary, we can assume that  MSCQ holds at every point $u\in \Gamma\cap U$ with modulus $\gamma>0.$
Take any $x\in U\cap \Gamma$ and $x^*\in \R^n.$
\par $(i)$ Suppose that $ v\in \R^n$  satisfying \eqref{msid}.  We first justify  $\Lambda(x, x^*; v)=\Lambda(x, x^*)$ by following the proof scheme of  \cite[Proposition 5.3]{GM15} with suitable modifications.
 If  either $\Lm(x, x^*)=\emptyset$ or  ${\rm span}\{\nabla q_i(x)|\; i\in E\cup I^+(x,x^*)\}=\{0\}$ then the assertion holds trivially. Assume now that  $\Lm(x, x^*)\neq\emptyset$ and ${\rm span}\{\nabla q_i(x)|\; i\in E\cup I^+(x,x^*)\}\not=\{0\}$. Let  $\widetilde E\subset E$ and $\tilde I\subset I^+(x,x^*)$ be such that  $\{\nabla q_j(x)|\; j\in B\}$ is a basis of ${\rm span}\{\nabla q_i(x)|\; i\in E\cup I^+(x,x^*)\}$ with $B:=\E\cup \tilde I\not=\emptyset$.
Consider the equations
 \begin{equation*}\label{imfu}
	Q_j(t, z):=q_j(x+tv+A^*z)=0,\; j\in B,
\end{equation*}
where  the rows of the $|B|\times n$ matrix $A$ are given by the gradients $\nabla q_j(x)$ for  $j\in B$, and $A^*$ is the transpose of $A$.
 Let $Q: \mathbb{R}\times \mathbb{R}^{|B|}\to \mathbb{R}^{|B|} $ be the function defined by $$Q(t,z)=\big(Q_j(t, z)\big)_{j\in B} \quad \text{ for all } t\in\mathbb{R}\  \mbox{and}\  z\in \mathbb{R}^{|B|}.$$ Then  $Q(0,0)=0$ and $\nabla_z Q(0, 0) =AA^*.$ Moreover,   the  $|B|\times |B|$ matrix  $AA^*$  is  invertible  since the rows of $A$ are linearly independent. Thus, by  the classical implicit function theorem \cite[p.~5]{DR14}, there exist  $\bar t>0$ and a $C^1$-smooth function $z: (-\bar t, \bar t)\to \R^{|B|}$ such that  $z(0)=0$ and
 \begin{equation}\label{imfu1}
Q_j\big(t, z(t)\big)=q_j\big(x+tv+A^*z(t)\big)=0 \quad  \text{ for all }\  t\in (-\bar t, \bar t)\ \mbox{and} \  j\in B,
\end{equation}
 Put  $\tilde{x}(t) := x + tv + A^*z(t)$ for $t\in  (-\bar t, \bar t).$ By \eqref{msid} and  \eqref{imfu1},  we get
$$0=\frac{d}{dt}q_j(\tilde{x}(t))\Big|_{t=0}=\Big\la\nabla q_j(x), v+A^*\frac{dz}{dt}(0)\Big\ra=\Big\la\nabla q_j(x), A^*\frac{dz}{dt}(0)\Big\ra,  j\in B. $$
 This shows that $AA^*\frac{dz}{dt}(0)=0$, and hence $\frac{dz}{dt}(0)=0$ since $AA^*$ is invertible. Therefore,
 \begin{equation}\label{imfu2}
\tilde{x}(0)=x,\; \frac{d\tilde{x}}{dt}(0)= v+ A^*\frac{dz}{dt}(0)=v, \text{ and } q_j\big(\tilde{x}(t)\big)=0\quad \mbox{for all}\  t\in (-\bar t, \bar t)\ \mbox{and} \  j\in B.
\end{equation}
Since $U$ is open, $x\in U$, $\lim\limits_{t\rightarrow 0}\tilde x(t)=x$, $\{\nabla q_j(x)|\; j\in B\}$ is a basis of ${\rm span}\{\nabla q_i(x)|\; i\in E\cup I^+(x,x^*)\},$ and $\{\nabla q_i(u)|\; i\in E\cup I^+(x,x^*)\}$ has the same rank for every $u\in U$ due to $I^+(x,x^*)\subset {I}(x)\subset {I}(\bar x)$, choosing $\bar t>0$  smaller  if necessary,   $\{\nabla q_j\big(\tilde{x}(t)\big)|\; j\in B\}$ is a basis of ${\rm span}\{\nabla q_i\big(\tilde{x}(t)\big)|\; i\in E\cup I^+(x,x^*)\}$ whenever $t\in (-\bar t, \bar t).$
Thus,   for every  $i\in \big(E\cup I^+(x,x^*)\big)\setminus B$ and  $t\in (-\bar t, \bar t)$ there exists $\alpha_{ij}(t)\in \R$ for $j\in B$ such that
 $$ \nabla q_i(\tilde{x}(t)) = \sum\limits_{j\in B}\alpha_{ij}(t)\nabla q_j(\tilde{x}(t)).$$
 This  together with \eqref{imfu2} shows  that
\begin{eqnarray*}
	\frac{d}{dt}q_i\big(\tilde{x}(t)\big)=\Big\la\nabla q_i(\tilde{x}(t)), \frac{d}{dt}\tilde{x}(t)\Big\ra	= \sum\limits_{j\in B}\alpha_{ij}(t)\Big\la\nabla q_j(\tilde{x}(t)), \frac{d}{dt}\tilde{x}(t)\Big\ra=0,
\end{eqnarray*}
for all $ t\in (-\bar t, \bar t)$  and  $i\in \big(E\cup I^+(x,x^*)\big)\setminus B.$ Thus, noting that $q_i(\tilde{x}(0))=0$ for all  $i\in \big(E\cup I^+(x,x^*)\big)\setminus B,$ we have
$q_i\big(\tilde{x}(t)\big)=0$ for all  $ t\in (-\bar t, \bar t)$  and  $i\in \big(E\cup I^+(x,x^*)\big)\setminus B.$
Combining this with \eqref{imfu2} tells us that
\begin{equation}\label{eqtv1}
q_i(\tilde{x}(t))=0\quad  \mbox{for all}\  i\in E\cup I^+(x,x^*)\ \mbox{and}\   t\in (-\bar t, \bar t).
\end{equation}
Take any $\lambda^{(1)},\lambda^{(2)}\in \Lambda(x,x^*).$ Since $I^+(\lambda^{(1)})\cup  I^+(\lambda^{(2)})\subset I^+(x,x^*)$,  we have  $\lambda_i^{(1)}=\lambda_i^{(2)}=0$ for all $i\in I\setminus I^+(x,x^*).$  This along with \eqref{eqtv1} gives us that
$\left\la \lambda^{(1)}-\lambda^{(2)}, q\big(\tilde{x}(t)\big)\right\ra = 0$
for all  $t\in (-\bar t, \bar t).$
On the other hand, we have
$$\begin{array}{rl}
&\lim\limits_{t\to 0}\dfrac{\left\la\lambda^{(1)}-\lambda^{(2)}, \nabla q(x)\big(\tilde{x}(t)-x\big)\right\ra}{t^2}=\lim\limits_{t\to 0}\dfrac{\left\la\nabla q(x)^*(\lambda^{(1)}-\lambda^{(2)}), \tilde{x}(t)-x\right\ra}{t^2}=0,\\ \cr
&\lim\limits_{t\to 0}\dfrac{\tilde{x}(t)-x}{t}=\lim\limits_{t\to 0}\dfrac{d}{dt}\tilde{x}(t)=v\quad\mbox{and}\quad  \lim\limits_{t\to 0}\dfrac{o(\|\tilde{x}(t)-x\|^2)}{t^2}=0.
\end{array}$$
Therefore,
\begin{eqnarray*}
	0&=&\lim\limits_{t\to 0}\frac{\left\la\lambda^{(1)}-\lambda^{(2)}, q\big(\tilde{x}(t)\big)-q\big(\tilde{x}(0)\big) \right\ra}{t^2}=\lim\limits_{t\to 0}\frac{\left\la\lambda^{(1)}-\lambda^{(2)}, q\big(\tilde{x}(t)\big)-q(x) \right\ra}{t^2}\\\cr
	&=&\lim\limits_{t\to 0} \scalebox{1.25}{$\frac{\la\lambda^{(1)}-\lambda^{(2)}, \nabla q(x) (\tilde{x}(t)-x)\ra +\frac{1}{2}\la\tilde{x}(t)-x, \nabla^2\la\lambda^{(1)}-\lambda^{(2)}, q\ra(x)(\tilde{x}(t)-x)\ra+o(\|\tilde{x}(t)-x\|^2)}{t^2}$}\\\cr
	&=& \frac{1}{2}\la v, \nabla^2\la\lambda^{(1)}-\lambda^{(2)}, q\ra(x)v\ra.
\end{eqnarray*}
 This shows that the form $\la v, \nabla^2\la \cdot, q\ra (x)v\ra $ is constant on $\Lambda(x, x^*)$, and   $\Lambda(x, x^*; v)=\Lambda(x, x^*)$ for every  $v\in \R^n$ satisfying (\ref{msid}).
 \par $(ii)$  Since  MSCQ holds at  $x,$ by  Lemma~\ref{cc4},  we see that $v$ fulfills \eqref{msid} whenever $v\in K(x, x^*)$. By $(i)$,
$\Lambda(x, x^*; v)=\Lambda(x, x^*)$   for every   $v\in K(x, x^*)$.
Therefore, by \cite[Theorem 3.5]{CH17}, we have
  \begin{equation}\label{tc}
  \begin{array}{rl}
T_{{\rm gph}N_{\Gamma}}(x, x^*)&=\big\{(w, z)\in\R^{2n}\big|\; \exists \lambda\in \Lambda(x, x^*; w) :  z\in \nabla^2\la \lambda, q\ra (x)w+ N_{K(x,x^*)}(w)\big\}\\ \cr
&=\big\{(w, z)\in\R^{2n}\big|\; \exists \lambda\in \Lambda(x, x^*) :  z\in \nabla^2\la \lambda, q\ra (x)w+ N_{K(x,x^*)}(w)\big\}.
\end{array}
\end{equation}
Let us now fix an arbitrary multiplier $\lambda\in \Lambda(x, x^*).$
By \eqref{tc}, we see that
$$\big\{(w, z)\in\R^{2n}\big|\; z\in \nabla^2\la \lambda, q\ra (x)w+ N_{K(x,x^*)}(w)\big\}\subset T_{{\rm gph}N_{\Gamma}}(x, x^*).$$
Conversely, take any $(w,z)\in T_{{\rm gph}\partial\delta_{\Gamma}}(x, x^*).$ Using \eqref{tc} one can find $\tilde\lambda\in \Lambda(x, x^*)$ such that
 $$z\in \nabla^2\la \tilde\lambda, q\ra (x)w+ N_{K(x,x^*)}(w).$$
 Since $\Lambda(x, x^*; v)=\Lambda(x, x^*)$  for every   $v\in K(x, x^*)$,
by Lemma~\ref{GO16lm4}, we have
$$\nabla^2\la \tilde\lambda-\lambda,q\ra(x)\in {\rm span}\{\nabla q_i(x)^*\ |\ i\in E\cup I^+(x,x^*)\}.$$
This along with ${\rm span}\{\nabla q_i(x)^*\ |\ i\in E\cup I^+(x,x^*)\} + N_{K(x,x^*)}(w)= N_{K(x,x^*)}(w)$ tells us  that
$$\begin{array}{rl}z&\in \nabla^2\la \tilde\lambda, q\ra (x)w+ N_{K(x,x^*)}(w)\\ \cr
&\subset \nabla^2\la \lambda, q\ra (x)w+ {\rm span}\{\nabla q_i(x)^*\ |\ i\in E\cup I^+(x,x^*)\} + N_{K(x,x^*)}(w)\\ \cr
&=\nabla^2\la \lambda, q\ra (x)w+N_{K(x,x^*)}(w),
\end{array}$$
which justifies that
$$T_{{\rm gph}N_{\Gamma}}(x, x^*)\subset\big\{(w, z)\in\R^{2n}\big|\; z\in \nabla^2\la \lambda, q\ra (x)w+ N_{K(x,x^*)}(w)\big\}.$$
Therefore,
$$T_{{\rm gph}N_{\Gamma}}(x, x^*)=\big\{(w, z)\in\R^{2n}\big|\; z\in \nabla^2\la \lambda, q\ra (x)w+ N_{K(x,x^*)}(w)\big\},$$
or equivalently,
$$
DN_\Gamma(x,x^*)(w)=\nabla^2\la \lambda, q\ra (x)w+ N_{K(x,x^*)}(w)\quad \mbox{for all}\ w\in \R^n,
$$
with an arbitrary multiplier $\lambda\in \Lm(x,x^*)$.
\par $(iii)$ Suppose further that $ \Lm(x,x^*)\not=\emptyset.$ Let $w\in \R^n$ and $z\in DN_\Gamma(x,x^*)(w).$ Since MSCQ holds at $x$ with modulus $\gamma$, by \cite[Lemma 3.2]{CH17}, we see that  $\Lambda(x, x^*)\cap\gamma\|x^*\|\B_{\R^\ell}\not=\emptyset.$ This allows us to pick
$\lambda\in \Lambda(x, x^*)\cap\gamma\|x^*\|\B_{\R^\ell}.$ Therefore, by $(ii)$, we have
$$z\in \nabla^2\la \lambda, q\ra (x)w+ N_{K(x,x^*)}(w),$$
which verifies the validity of $(iii)$.
$\hfill\Box$

 \section{Second-order Characterizations of Tilt Stability for NLPs under RCRCQ}

This section is devoted to the study of tilt stability of the classical nonlinear program  \eqref{tiltMP1}, which is formulated as follows:
   \begin{equation*} \begin{array}{rl} &\mbox{minimize} \quad \varphi(x)\quad \mbox{subject to} \quad q_i(x)=0,\;i\in E \text{ and } q_i(x)\leq 0,\ i\in I,
 	\end{array}
 \end{equation*}
 where $\varphi:\R^n\rightarrow\R$ and $q_i: \R^n\to \R$ is a twice continuously differentiable function for every $i\in E\cup I$. Here  $E := \{i\in \N\ |\ 1\leq i\leq {\ell}_1\}$ and $I := \{i\in \N\ |\ {\ell}_1 + 1\leq i\leq {\ell}_1 + {\ell}_2\}$, with $\ell_1,\ell_2\in \N:=\{0,1,2,...\}$ and $E\cup I\not=\emptyset$.

The above  nonlinear program  can be written as the unconstrained optimization problem:
 \begin{equation*}\label{NP2} \mbox{minimize} \quad f(x):=\varphi(x)+\delta_\Gamma(x)\; \text{ on } \mathbb{R}^n,
 \end{equation*}
 where $\Gamma$  is  the feasible set of   \eqref{tiltMP1}  given  by \eqref{CS},  and $\delta_\Gamma$ is the indicator function of the set $\Gamma\subset \R^n$. 

Recall \cite{M24} that  a point $\ox\in \Gamma$ is  said to be  a  {\it tilt-stable local minimizer} of   problem \eqref{tiltMP1} with
 modulus $\kk>0$  if  it is a  tilt-stable local minimizer of the function $f$ defined above, that is,
 there exists $\gg>0$ such that the argmin solution mapping
 \[
  M_\gg(v):=\disp{\rm argmin}\,\Big\{\varphi(x)-\la v,x\ra|\;  x\in \Gamma\cap \B_\gg(\ox)\Big\}
  \]
 is single-valued and Lipschitz continuous with constant $\kk>0$ on some neighborhood of $0\in \R^n$ with $M_\gg(0)=\ox$.
 One says that $\ox$ is a  {\it tilt-stable local minimizer} of  problem \eqref{tiltMP1} if it is a  tilt-stable local minimizer of  problem \eqref{tiltMP1} with some modulus $\kappa>0.$
 The number  ${\rm tilt}(\varphi,q,\ox):={\rm tilt}(f,\ox)$ is called the {\it exact bound of tilt stability} of problem  \eqref{tiltMP1} at  $\ox$.

  For $x\in \Gamma$ and $v\in \R^n,$ from \eqref{som} it follows  that
 $$\begin{array}{rl}\Lm\big(x,v-\nabla \varphi (x)\big)&=\left\{\lambda\in N_\Theta\big(q(x)\big)\ |\ \nabla q(x)^*\lambda=v-\nabla \varphi (x)\right\}\\
 &=\left\{\lambda\in \R^{\ell_1}\times \R_+^{\ell_2}\ | \nabla_x\mathcal{L}(x, \lambda)=v\ \mbox{and}\    \lambda_iq_i(x)=0\quad \mbox{for all} \ i\in I\right\}.
 \end{array}
 $$
In particular, we see that $\Lm\big(\bar x,-\nabla \varphi (\bar x)\big)=\Lambda(\bar x)$ for $\bar x\in \Gamma$.

\begin{lemma}{\rm (\cite[Lemma 3.2]{AMS05}).} \label{lem1} Let  $f_0, f_1,...,f_q : \R^n \to \R$  be
 functions that are continuously differentiable on some neighborhood of  $\bar x\in \R^n$ with  $\{\nabla f_1(\ox), . . . ,\nabla f_q (\ox)\}$  being  linearly
independent. Suppose further that   $\nabla f_0 (x) \in
{\rm span}\{\nabla f_1(x), . . . ,\nabla f_q (x)\}$  for all $x\in\R^n$ near  $\ox$.  Then there exist open convex neighborhoods $X$ of $\ox$ and $Y$ of $\big( f_1( \ox), . . . , f_q ( \ox)\big)$ along with  a continuously differentiable function $g_0 : Y \to \R$  such that $\big( f_1(x), . . . , f_q (x)\big)\in Y$ and  $f_0(x) = g_0\big( f_1(x), . . . , f_q (x)\big)$ for every  $x \in X.$  Moreover, $\nabla f_0 ( \ox) =\sum\limits_{i=1}^q\lambda_i\nabla f_i ( \ox),$ where
$\lambda_i:=\frac{\partial g_0}{\partial y_i}\big( f_1( \ox), . . . , f_q ( \ox)\big)$ for $i=1,2,...,q.$
\end{lemma}

The next result is used to derive  new characterizations of tilt stability in what follows.

\begin{lemma} \label{prop1shulu}	Let  RCRCQ hold  at $\bar x\in \Gamma$, and    $\bar\lambda, \;\tilde{\lambda}\in \Lambda\big(\ox, -\nabla \varphi(\ox)\big)$ with $I^+(\tilde{\lambda})=I^+\big(\bar x,-\nabla \varphi(\ox)\big)$, where  $I^+\big(\bar x, -\nabla \varphi(\ox)\big)$ is defined as in \eqref{icong} with $x=\bar x$ and $x^*=-\nabla\varphi (\ox)$ and $I^+(\lambda)$ is  taken from \eqref{sci}. Assume further that there exist sequences $x_k\st{\Gamma}{\to} \ox,\; v_k\to 0$  as $k\to \infty$ and $\lambda^k\in \Lambda\big(x_k, v_k-\nabla \varphi(x_k)\big)$ such that $\lambda^k\to \bar\lambda$ as $k\to \infty$. Then we have the equality
\begin{equation}\label{equa7}	\sum\limits_{i\in I^+\big(\bar x,-\nabla\varphi(\bar x)\big)\setminus I^+(\bar\lambda) }\tilde{\lambda}_i\nabla q_i(\ox) + \sum\limits_{i\in E\cup I^+(\bar\lambda) }(\tilde{\lambda}_i-\bar\lambda_i)\nabla q_i(\ox) = 0,
\end{equation}
and $I^+\big(\bar x,-\nabla\varphi(\bar x)\big)\subset {I}(x_k)$ for all $k$ sufficiently large.
\end{lemma}
\noindent{\it Proof.} Since  $\bar\lambda, \;\tilde{\lambda}\in \Lambda\big(\ox, -\nabla \varphi(\ox)\big)$ with $I^+(\tilde{\lambda})=I^+\big(\bar x,-\nabla\varphi(\bar x)\big)$, we see  that
$$I^+({\bar\lambda})\subset I^+\big(\bar x,-\nabla\varphi(\bar x)\big)$$
 and
 $$\sum\limits_{i\in E\cup I^+({\bar\lambda})}\bar\lambda_i\nabla q_i(\ox)=\nabla q(\ox)^*\bar\lambda=-\nabla\varphi (\ox)=\nabla q(\ox)^*\tilde\lambda=\sum\limits_{i\in E\cup I^+\big(\bar x,-\nabla\varphi(\bar x)\big)}\tilde\lambda_i\nabla q_i(\ox).$$
 Thus, we get  \eqref{equa7}.
Pick $\E\subset E,$   $I_0\subset I^+(\bar\lambda),$ and $J_0\subset I^+\big(\bar x,-\nabla\varphi(\bar x)\big)\setminus I^+(\bar\lambda)$   such that $\big\{\nabla q_i(\ox)\ |\ i\in \E\cup I_0\big\}$ and $\{\nabla q_i(\ox)\ |\  i\in L\}$ are  bases of ${\rm span}\big\{\nabla q_i(\ox)\ |\ i\in E\cup I^+(\bar\lambda)\big\}$ and  ${\rm span}\{\nabla q_i(\ox)\ |\ i\in E\cup I^+\big(\bar x,-\nabla\varphi(\bar x)\big)\}$, respectively, where $L:=\E\cup I_0\cup J_0$. Then,  by  \eqref{equa7}, for $i\in \tilde E\cup I_0$  one can find  $\mu_i\in \R$  such that
\begin{equation}\label{equa6}
	\sum\limits_{i\in J_0}\tilde{\lambda}_{i}\nabla q_{i}(\ox) +\sum\limits_{i\in \E\cup I_0}\mu_i\nabla q_{i}(\ox) + \sum\limits_{j\in K}\tilde{\lambda}_{j}\nabla q_{j}(\ox)= 0, \; \tilde{\lambda}_i>0  \mbox{ as }i\in J_0\cup K,
\end{equation}
where $K:=I^+\big(\bar x,-\nabla\varphi(\bar x)\big)\setminus \big(J_0\cup I^+(\bar\lambda)\big).$    Since  RCRCQ holds at $\bar x,$ $K\subset E\cup I^+\big(\bar x,-\nabla\varphi(\bar x)\big),$   and $\{\nabla q_i(\ox)\ |\  i\in L\}$  is a basis of ${\rm span}\left\{\nabla q_i(\ox)\ |\ i\in E\cup I^+\big(\bar x,-\nabla\varphi(\bar x)\big)\right\},$ it follows that $\{\nabla q_i(x)\ |\  i\in L\}$ is linearly independent, and
$\nabla q_j(x)\in {\rm span}\{\nabla q_i(x)\ |\  i\in L\}$ for all $x$ near $\ox$ and for all $j\in K.$ By Lemma~\ref{lem1}, there exist open convex neighborhoods $X$ of $\ox$ and  $Y$ of $q_{L}(\ox)\in \R^{|L|}$ along with  continuously differentiable functions $g_j: Y\to \R,$ $j\in K,$ such that $q_{L}(x)\in Y$ and
\begin{equation} \label{crt1}
	q_{j}(x)=g_j\big(q_{L}(x)\big) \quad \mbox{for all}\  x\in X,
\end{equation}
where $q_{L}: X\to \R^{|L|}$ is a mapping with component functions being $q_i(x)$ for $i\in L$. Therefore,
$$ \nabla q_{j}(x)=\sum\limits_{i\in L}\frac{\partial g_j}{\partial y_i}\big(q_{L}(x)\big)\nabla q_i(x)\quad \mbox{for all}\ x\in X \mbox{ and } j\in K.$$
In particular,
$$ \nabla q_{j}(\ox)=\sum\limits_{i\in L}\frac{\partial g_j}{\partial y_i}\big(q_{L}(\ox)\big)\nabla q_i(\ox)=\sum\limits_{i\in \E\cup I_0\cup J_0}\frac{\partial g_j}{\partial y_i}\big(q_{L}(\ox)\big)\nabla q_i(\ox)\quad \mbox{for all}\ j\in K.$$
This along with \eqref{equa6}  gives us the equality
	$$\sum\limits_{i\in J_0}\left(\tilde{\lambda}_{i}+\sum\limits_{j\in K}\tilde\lambda_j\frac{\partial g_j}{\partial y_i}\big(q_{L}(\ox)\big)\right)\nabla q_{i}(\ox) +\sum\limits_{i\in \E\cup I_0}\left(\mu_i+\sum\limits_{j\in K}\tilde\lambda_j\frac{\partial g_j}{\partial y_i}\big(q_{L}(\ox)\big)\right)\nabla q_{i}(\ox) = 0.$$
Hence, noting that $\{\nabla q_i(\bar x)\ |\ i\in \E\cup I_0\cup J_0\}$ is linearly independent, we have
\begin{equation} \label{Sh-eq1}\sum\limits_{j\in K}\tilde\lambda_j\frac{\partial g_j}{\partial y_i}\big(q_{L}(\ox)\big)=-\tilde\lambda_{i}<0\quad \mbox{for all}\  i\in J_0.\end{equation}
Following the proof scheme of Lu \cite[Proposition 1]{ShuLu}, we next show that $I^+\big(\bar x,-\nabla\varphi(\bar x)\big)\subset {I}(x_k)$ for all $k$ sufficiently large.
Put $g(y)=\sum\limits_{j\in K}\tilde\lambda_j g_j(y)$ for $y\in Y.$ By \eqref{Sh-eq1}, we get
$$\frac{\partial g}{\partial y_i}\big(q_{L}(\ox)\big)=\sum\limits_{j\in K}\tilde\lambda_j\frac{\partial g_j}{\partial y_i}\big(q_{L}(\ox)\big)<0\quad \mbox{for}\ i\in J_0.$$
Thus, shrinking $X$ and $Y$  if necessary, we may assume that $\frac{\partial g}{\partial y_i}\big(y\big)<0\quad \mbox{for all}\ y\in Y\ \mbox{and}\ i\in J_0.$
Note that $q_i(\ox)=0$ for every $i\in E\cup I^+\big(\bar x,-\nabla\varphi(\bar x)\big)$ and  $K\cup L \subset E\cup I^+\big(\bar x,-\nabla\varphi(\bar x)\big).$  So, by~\eqref{crt1},
we have $0=q_L(\ox)\in Y$ and
$$g(0)=g\big(q_L(\ox)\big)=\sum\limits_{j\in K}\tilde\lambda_j g_j\big(q_L(\ox)\big)=\sum\limits_{j\in K}\tilde\lambda_j q_j(\ox)=0.$$
Since $\lim\limits_{k\to \infty}x_k=\bar x\in X,$ $X$ is open,  $\lim\limits_{k\to \infty}\lambda^k=\bar\lambda$  and $I_0\subset I^+(\bar\lambda)$, we see that  $I_0\subset I^+(\lambda^k)$ and  $x_k\in X,$ and thus
 \begin{equation}\label{eqtv2} q_L(x_k)\in Y\ \mbox{and}\   q_i(x_k)=0,\end{equation}
for all $i\in  \E\cup I_0=L\backslash J_0$ and for all $k$ sufficiently large.
 We claim that \eqref{eqtv2} also holds for every $i\in J_0$ and for all $k$ sufficiently large.
 Suppose on the contrary that the claim fails. Then
 there exist $i_0\in J_0$ and  $k_0$ sufficiently large such that  $q_L(x_{k_0})\in Y$ and $q_{i_0}(x_{k_0})<0$.
 Thus, by the classical mean value theorem, we see that
\begin{eqnarray}
	g\big(q_{L}(x_{k_0})\big) &=& g\big(q_{L}(x_{k_0})\big) - g(0) = \big\la\nabla g(y^*_{k_0}), q_{L}(x_{k_0})\big\ra\notag\\
	&=&  \sum\limits_{i\in J_0}\frac{\partial g}{\partial y_i}(y^*_{k_0})q_{i}(x_{k_0})+\sum\limits_{i\in L\backslash J_0}\frac{\partial g}{\partial y_i}(y^*_{k_0})q_{i}(x_{k_0}), \label{16624}
\end{eqnarray}
for some $y^*_{k_0} $ lying between $0=q_L(\ox)$ and $q_{L}(x_{k_0})$.  Since $Y$ is convex, $0=q_L(\ox)\in Y,$  $q_{L}(x_{k_0})\in Y$, and  $y^*_{k_0} $ belongs to the interval  $[ q_L(\ox), q_{L}(x_{k_0})]$, we get   $y^*_{k_0}\in Y,$ and thus $\frac{\partial g}{\partial y_i}\big(y^*_{k_0}\big)<0$ for all  $i\in J_0.$
 Noting that $\tilde\lambda_i>0$ and  $q_{i}(x_{k_0})\leq 0$ for every $i\in J_0,$  $q_{i_0}(x_{k_0})<0$, $i_0\in J_0,$   and $q_{i}(x_{k_0})=0$ for every $i\in L\backslash J_0$, from \eqref{16624} it follows that  $g\big(q_{L}(x_{k_0})\big)>0.$
  On the other hand, by the definition of $g$ along with  \eqref{crt1}, we have
\begin{equation*}\label{166a24}
	g\big(q_{L}(x_{k_0})\big) = \sum\limits_{j\in K}\tilde\lambda_j g_j\big(q_{L}(x_{k_0})\big)= \sum\limits_{j\in K}\tilde\lambda_jq_j(x_{k_0})\leq 0.
\end{equation*}
 Therefore,  we  get a contradiction, which shows that $q_{i}(x_{k})=0$ for $i\in J_0$ and $k$ sufficiently large. For every $k$ sufficiently large, combining the latter with \eqref{eqtv2} tells us that $q_{i}(x_{k})=0=q_{i}(\bar x)$ for $i\in L,$ and hence
 $$q_{j}(x_{k})=g_j\big(q_L(x_k)\big)=g_j\big(q_L(\bar x)\big)=q_j(\bar x)=0\quad \mbox{for all}\ j\in K.$$
So, noting that  $I^+\big(\bar x,-\nabla\varphi(\bar x)\big)=K\cup J_0\cup I^+(\bar \lambda)$ and  $I^+(\bar \lambda)\subset {I}(x_k)$ for every $k$ sufficiently large, we see that    $I^+\big(\bar x,-\nabla\varphi(\bar x)\big)\subset {I}(x_k)$ for every $k$ sufficiently large. $\hfill\Box$


The main result of this paper is the following theorem, which provides some pointbased second-order characterizations of tilt stability along with  a formula for calculating the exact bound of tilt stability for NLPs under RCRCQ.

\begin{theorem}\label{sochots}  Let $\ox\in \Gamma$ be a stationary point of  \eqref{tiltMP1} at which RCRCQ is fulfilled. Then the following assertions hold:
\par	$(i)$ Given $\kk>0$, the necessary and sufficient condition for $\ox$ to be  a tilt-stable local minimizer of~\eqref{tiltMP1} with any modulus $\kk'>\kk$ is that
\begin{equation}\label{soch1}\big\la w, \nabla_x^2\mathcal{L}(\ox, \lambda)w\big\ra  \ge\frac{1}{\kk}\|w\|^2\quad \text{whenever } \lambda\in \Lambda(\ox)\ \mbox{and}\  	\la\nabla q_i(\ox), w\ra =0\ \mbox{as}\  i\in E\cup I^+,	\\
\end{equation}
where $I^+:=I^+\big(\bar x, -\nabla\varphi(\bar x)\big)$ and $\Lambda(\ox)=\Lambda\big(\ox, -\nabla \varphi(\ox)\big)$.    Furthermore, condition \eqref{soch1}  amounts to  that there exists  $\lambda\in \Lambda(\ox)$ such that
\begin{equation}\label{soch2} \big\la w, \nabla_x^2\mathcal{L}(\ox, \lambda)w\big\ra  \ge\frac{1}{\kk}\|w\|^2 \quad \text{whenever } \  \la\nabla q_i(\ox), w\ra =0\ \mbox{as}\  i\in E\cup I^+.
\end{equation}
\par$(ii)$   The point  $\ox$ is a tilt-stable local minimizer of \eqref{tiltMP1}  if and only if for every nonzero vector $w\in \R^n$ with $\la\nabla q_i(\bar x), w\ra =0$ as $i\in E\cup I^+$ we have
\begin{equation}\label{soch3} \big\la w, \nabla_x^2\mathcal{L}(\ox, \lambda)w\big\ra  >0\quad \text{for all}\  \lambda\in \Lambda(\ox\big),\\
\end{equation}
which is equivalent to the  condition: there exists  $\lambda\in \Lambda(\ox)$ such that $\la w, \nabla_x^2\mathcal{L}(\ox, \lambda)w\ra>0$
for every nonzero vector $w\in \R^n$ with $\la\nabla q_i(\bar x), w\ra =0$ as $i\in E\cup I^+$.

\par$(iii)$  If $\bar x$ is a tilt-stable local minimizer of \eqref{tiltMP1}, then for any $\bar\lambda\in \Lambda(\bar x)$  the exact bound of tilt stability is  calculated by
\begin{eqnarray}\label{tiltEB}
	{\rm tilt}\, (\varphi, q,\ox)=\sup\limits_{w}\left\{\frac{\|w\|^2} {\la w, \nabla_x^2\mathcal{L}(\ox, \bar\lambda)w\ra}\Big|\;   	\big\la\nabla q_i(\bar x), w\big\ra =0, i\in E\cup I^+\right\},
\end{eqnarray}
with the convention that $0/0:=0$.
\end{theorem}
\noindent{\it Proof.} $(i)$ We begin with  justifying the sufficient condition of the first statement in $(i)$. To this end,  let  \eqref{soch1} be fulfilled. Take any   $\tilde{\lambda}\in \Lambda\big(\ox, -\nabla \varphi(\ox)\big)$ with $I^+(\tilde{\lambda})=I^+$ which can be done due to Lemma~\ref{cc4}. Suppose on the contrary that  $\bar x$ is not  a tilt-stable local minimizer of \eqref{tiltMP1} with some modulus $\kk'>\kk$.  By Theorem~$\ref{theo22}$, there exist sequences $x_k\to \ox,\; v_k\to 0$ as $k\to\infty$ and $(z_k, w_k)\in T_{{\rm gph}\partial f}(x_k, v_k)$ such that
\begin{equation} \label{inequa1}
	\la z_k, w_k\ra  < \frac{1}{\kk'}\|w_k\|^2\quad  \text{for  all } k.
\end{equation}
This implies that $w_k\neq 0$  for every $k.$  Thus, replacing $(z_k, w_k)$ by $\frac{1}{\|w_k\|}(z_k, w_k)$ and passing to a subsequence if necessary, we can assume that  $(z_k, w_k)\in T_{{\rm gph}\partial f}(x_k, v_k)$ together with    $\|w_k\|=1$ for all $k,$ and  $\lim\limits_{k\to\infty}w_k= w$ for some vector $w\in \R^n$.   For every $k,$ by  the definition of subgradient graphical derivative and using \cite[Propostion 1.76]{M24},  we  get the following relations:
$$z_k\in \partial^2_Df(x_k,v_k)(w_k)=\nabla^2\varphi(x_k)w_k+\partial^2_D\delta_\Gamma\big(x_k, v_k-\nabla\varphi(x_k)\big)(w_k),$$
where $f(x):=\varphi(x)+\delta_{\Gamma}(x).$
This shows that
$$z_k-\nabla^2\varphi(x_k)w_k\in \partial^2_D\delta_\Gamma\big(x_k, v_k-\nabla\varphi(x_k)\big)(w_k)=DN_\Gamma\big(x_k, v_k-\nabla\varphi(x_k)\big)(w_k).$$
Furthermore,  $\Lambda\big(x_k, v_k-\nabla\varphi(x_k)\big)\not=\emptyset$ for every $k$ sufficiently large.
Therefore, by Theorem~\ref{prop53}, there exists $\lambda^k\in \Lambda\big(x_k, v_k-\nabla\varphi(x_k)\big)\cap \gamma\|v_k-\nabla\varphi(x_k)\|\B_{\R^\ell}$ such that
$$z_k-\nabla^2\varphi(x_k)w_k\in \nabla^2\big\la\lambda^k,q\big\ra(x_k)(w_k)+N_{K\big(x_k,v_k-\nabla\varphi(x_k)\big)}(w_k),$$
 or equivalently,
$$\begin{array}{rl}z_k-\nabla^2\varphi(x_k)w_k- \nabla^2\langle \lambda^k,q\rangle(x_k)w_k&\in N_{K\big(x_k,v_k-\nabla\varphi(x_k)\big)}(w_k)\\ \cr
&=\left[K\big(x_k,v_k-\nabla\varphi(x_k)\big)\right]^*\cap\{w_k\}^\perp,\end{array}$$
for all   $k$ sufficiently large.
This implies that  the sequence $\{\lambda^k\}$  is bounded, and
\begin{equation}\begin{array}{rl}\label{inequa2}\la z_k,w_k\ra&=\la \nabla^2\varphi(x_k)w_k+ \nabla^2\langle \lambda^k,q\rangle(x_k)w_k,w_k\ra\\ \cr
&=\langle w_k,\nabla^2_x\mathcal{L}(x_k,\lambda^k)w_k\rangle\quad \mbox{for all}\ k \ \mbox{sufficiently large}.\end{array}\end{equation}
Passing to a subsequence if necessary gives us  $\lambda^k\to \bar{\lambda}$ as $k\to\infty$ for some $\bar\lambda\in \R^\ell.$  Thus, by  $\lambda^k\in\Lambda\big(x_k,v_k-\nabla\varphi(x_k)\big)$, we get
$$\bar{\lambda}\in N_{\Theta}\big(q(\ox)\big)\ \mbox{and}\ \nabla q(\ox)^*\bar{\lambda} =\lim\limits_{k\to\infty } \nabla q(x_k)^*\lambda^k=\lim\limits_{k\to\infty } \big(v_k-\nabla \varphi(x_k)\big)=-\nabla\varphi(\ox),$$
which  justifies that  $\bar{\lambda}\in \Lambda\big(\ox, -\nabla\varphi(\ox)\big)$.
 This along with Lemma \ref{prop1shulu} gives  us that  $I^+\subset {I}(x_k)$ for all $k$ sufficiently large.
 On the other hand, since $w_k\in K\big(x_k,v_k-\nabla\varphi(x_k)\big)$ and $\lambda^k\in \Lambda\big(x_k, v_k-\nabla\varphi(x_k)\big)$, by Lemma \ref{cc4}, we see that
 \[\la\nabla q_i(x_k), w_k\ra\begin{cases}
	=0 &\text{ if } i\in E\cup I^+(\lambda^k),\\
	\leq 0&\text{ if } i\in  {I}(x_k)\setminus I^+(\lambda^k),\\
\end{cases}\]
for all $k$ sufficiently large. Therefore,  letting $k\to \infty$,   we get
\begin{equation}\label{tvc1}	\la\nabla q_i(\ox), w\ra\begin{cases}
	=0 &\text{ if } i\in E\cup I^+(\bar\lambda),\\
	\leq 0&\text{ if } i\in  {I}(\ox)\setminus I^+(\bar\lambda).
\end{cases}\end{equation}
Moreover, according to \eqref{equa7},
\begin{equation*}	\sum\limits_{i\in I^+\setminus I^+(\bar\lambda) }\tilde{\lambda}_i\nabla q_i(\ox) + \sum\limits_{i\in E\cup I^+(\bar\lambda) }(\tilde{\lambda}_i-\bar\lambda_i)\nabla q_i(\ox) = 0,
\end{equation*}
which along with  \eqref{tvc1}, $I^+\subset I^+(\bar x)$ and  $\tilde{\lambda}_i>0$ as $i\in I^+\setminus I^+(\bar\lambda)$ gives us that
$$ \la\nabla q_i(\ox), w\ra = 0 \text{ for } i\in I^+\setminus I^+(\bar\lambda)$$
So, we have
$$ \la\nabla q_i(\ox), w\ra = 0 \text{ for every}\  i\in E\cup I^+.$$
Combining \eqref{inequa1} and \eqref{inequa2} gives us that
$$\la w_k, \nabla^2_x\mathcal{L}(x_k,\lambda^k)w_k\ra = \la z_k, w_k\ra < \frac{1}{\kk'}\|w_k\|^2 \quad  \text{for all}\  \ k. $$
Passing to the limits as $k\to \infty$, we obtain
\begin{equation*}\label{inequa3}
	\la w, \nabla^2_x\mathcal{L}(\ox,\bar{\lambda})w\ra \leq  \frac{1}{\kk'}\|w\|^2 < \frac{1}{\kk}\|w\|^2,
\end{equation*}
where  $\bar{\lambda}\in \Lambda\big(\ox):=\Lambda\big(\ox, -\nabla\varphi(\ox)\big)$ and  $ \la\nabla q_i(\ox), w\ra = 0$  for every  $i\in E\cup I^+.$
This contradicts~\eqref{soch1}. Therefore, the point $\bar x$ is a tilt-stable local minimizer of~\eqref{tiltMP1} with any modulus $\kk'>\kk$.

We next prove the necessary condition. Suppose that the point $\bar x$ is a tilt-stable local minimizer of~\eqref{tiltMP1} with any modulus $\kk'>\kk$.  We need to show that \eqref{soch1} holds.
Suppose on the contrary that $\ox$ is a tilt-stable local minimizer of \eqref{tiltMP1}
with modulus $\kk'$ for every $\kk'>\kk$, but there are vectors $\lambda\in \Lambda\big(\ox, -\nabla \varphi(\ox)\big)$ and $w\in\R^n\setminus\{0\}$ such that $\la\nabla q_i(\ox), w\ra =0$  for every  $i\in E\cup I^+$ and
\begin{equation*} \la w, \nabla_x^2\mathcal{L}(\ox, \lambda)w\ra  <\frac{1}{\kk}\|w\|^2.
\end{equation*}
Then we can pick  $\kk'>\kk$ near $\kk$ such that
\begin{equation}\label{soch4} \la w, \nabla_x^2\mathcal{L}(\ox, \lambda)w\ra  <\frac{1}{\kk'}\|w\|^2.
\end{equation}
Since $\lambda\in \Lambda\big(\ox, -\nabla \varphi(\ox)\big)$ and  $\la\nabla q_i(\ox), w\ra =0$  for all   $i\in E\cup I^+$, by Lemma \ref{cc4}, we see that $w\in K\big(\bar x,-\nabla\varphi(\bar x)\big),$ and thus $0\in N_{ K\big(\bar x,-\nabla\varphi(\bar x)\big)}(w).$  Thus, by Theorem~\ref{prop53} and the sum rule for graphical derivative, we get
$$\begin{array}{rl} \nabla_x^2\mathcal{L}_x(\ox, \lambda)w&=\nabla^2\varphi(\bar x)w+\nabla^2\la \lambda,q\ra(\bar x)\in \nabla^2\varphi(\bar x)w+\nabla^2 \la \lambda,q\ra(\bar x)+N_{ K\big(\bar x,-\nabla\varphi(\bar x)\big)}(w)\\ \cr
&=\nabla^2\varphi(\bar x)w+DN_\Gamma\big(\bar x,-\nabla\varphi(\bar x)\big)(w)=\partial^2_Df\big(\bar x, 0)(w),
\end{array}$$
where $f(x):=\varphi(x)+\delta_{\Gamma}(x).$
Combining the latter with the assumption that $\ox$ is a tilt-stable local minimizer of \eqref{tiltMP1}
with modulus $\kk'$, by Theorem~\ref{theo22}, we have
$$\la w, \nabla_x^2\mathcal{L}(\ox, \lambda)w\ra  \ge \frac{1}{\kk'}\|w\|^2 ,$$
which contradicts \eqref{soch4}. This proves that  \eqref{soch1} holds.

Finally, we show the validity of the second statement in $(i)$.  From Theorem~\ref{prop53} it follows that
 $$\la w, \nabla_x^2\mathcal{L}(\ox, \lambda^1)w\ra=\la w, \nabla_x^2\mathcal{L}(\ox, \lambda^2)w\ra,$$
 whenever  $\lambda^1,\lambda^2\in \Lambda(\bar x)$  and    $\la\nabla q_i(\ox), w\ra =0$  as  $i\in E\cup I^+.$  This along with  $\Lambda(\bar x)\not=\emptyset$ shows that
\eqref{soch1} amounts to the existence of  $\lambda\in \Lambda\big(\ox\big)$  satisfying \eqref{soch2}.

$(ii)$ Suppose that $\ox$ is a tilt-stable local minimizer of \eqref{tiltMP1}. Then there exists $\kappa>0$ such that  $\ox$ is a tilt-stable local minimizer of \eqref{tiltMP1} with modulus $\kappa.$  This implies that  $\ox$ is a tilt-stable local minimizer of \eqref{tiltMP1} with any modulus $\kappa'>\kappa$. Thus, by $(i),$ we get \eqref{soch1}, which guarantees the validity of \eqref{soch3} for every nonzero vector $w\in \R^n$ with $\la\nabla q_i(\bar x), w\ra =0$ as $i\in E\cup I^+$.
Conversely, let \eqref{soch3} hold for every nonzero vector $w\in \R^n$ with $\la\nabla q_i(\bar x), w\ra =0$ as $i\in E\cup I^+$.
Let us consider the following two cases.
\par {\it Case 1:} $w=0$ whenever $\la\nabla q_i(\bar x), w\ra =0$ as $i\in E\cup I^+$. Then we see that  \eqref{soch1} holds for every $\kappa>0$. Hence, by $(i)$, the point $\bar x$ is a tilt-stable local minimizer of \eqref{tiltMP1}.
\par {\it Case 2:} There exists  a nonzero vector $w\in \R^n$ with $\la\nabla q_i(\bar x), w\ra =0$ as $i\in E\cup I^+$.  Since $\bar x$ is a stationary point of \eqref{tiltMP1}, one can find   $\bar\lambda\in \Lambda(\bar x)$.  Then from  \eqref{soch3} it follows that
\begin{equation*}\label{rewsoch3} \la w, \nabla_x^2\mathcal{L}(\ox, \bar{\lambda})w\ra  \geq \mu\|w\|^2 \text{ whenever } \la\nabla q_i(\ox), w\ra =0 \text{ as } i\in E\cup I^+	,\\
\end{equation*}
where  $\mu:=\inf\big\{\la w, \nabla_x^2\mathcal{L}(\ox, \bar{\lambda})w\ra\big|\;\|w\|=1, \la\nabla q_i(\ox), w\ra =0, i\in E\cup I^+ \big\} $ is a positive real number.
Therefore, by $(i)$, the point $\bar x$ is a tilt-stable local minimizer of \eqref{tiltMP1}.

$(iii)$ Suppose $\bar x$ is a  tilt-stable local minimizer of \eqref{tiltMP1}. Take any $\kappa> {\rm tilt}\, (\varphi, q,\ox).$ Then $\bar x$ is a  tilt-stable local minimizer of \eqref{tiltMP1} with any modulus $\kappa'>\kappa.$ Thus, for any $\bar\lambda\in \Lambda(\bar x)$, by $(i)$, we have
$$\big\la w, \nabla_x^2\mathcal{L}(\ox, \bar\lambda)w\big\ra  \ge\frac{1}{\kk}\|w\|^2 \quad \text{whenever } \  \la\nabla q_i(\ox), w\ra =0\ \mbox{as}\  i\in E\cup I^+,$$
and hence
$$\kappa\geq \sup\limits_{w}\Big\{\frac{\|w\|^2} {\la w, \nabla_x^2\mathcal{L}(\ox, \bar\lambda)w\ra}\Big|\;   	\big\la\nabla q_i(\bar x), w\big\ra =0, i\in E\cup I^+\Big\}.$$
Therefore, letting $\kappa\to {\rm tilt}\, (\varphi, q,\ox)$ with
$\kappa> {\rm tilt}\, (\varphi, q,\ox),$ we get
\begin{equation}\label{tiltmod1}
 {\rm tilt}\, (\varphi, q,\ox)\geq \sup\limits_{w}\Big\{\frac{\|w\|^2} {\la w, \nabla_x^2\mathcal{L}(\ox, \bar\lambda)w\ra}\Big|\;   	\big\la\nabla q_i(\bar x), w\big\ra =0, i\in E\cup I^+\Big\}.\end{equation}
  Conversely, take any $\bar\lambda\in \Lambda(\bar x)$. Since  $\bar x$ is a  tilt-stable local minimizer of \eqref{tiltMP1} with some modulus $\kappa_0>0$, by $(i)$, we have
 $$\la w, \nabla_x^2\mathcal{L}(\ox, \bar\lambda)w\ra\geq \frac{1}{\kappa_0}\|w\|^2 \quad \text{whenever } \  \la\nabla q_i(\ox), w\ra =0\ \mbox{as}\  i\in E\cup I^+.$$
 This implies that
\begin{equation}\label{tiltmod2}\la w, \nabla_x^2\mathcal{L}(\ox, \bar\lambda)w>0 \quad \text{for every  $w\in \R^n\backslash\{0\}$ with } \  \la\nabla q_i(\ox), w\ra =0\ \mbox{as}\  i\in E\cup I^+.\end{equation}
 Suppose on the contrary that
  $${\rm tilt}\, (\varphi, q,\ox)>\sup\limits_{w}\Big\{\frac{\|w\|^2} {\la w, \nabla_x^2\mathcal{L}(\ox, \bar\lambda)w\ra}\Big|\;   	\big\la\nabla q_i(\bar x), w\big\ra =0, i\in E\cup I^+\Big\}.$$
 Then we can find $\kappa',\kappa$ such that
$${\rm tilt}\, (\varphi, q,\ox)>\kappa'>\kappa>\sup\limits_{w}\Big\{\frac{\|w\|^2} {\la w, \nabla_x^2\mathcal{L}(\ox,\bar \lambda)w\ra}\Big|\;   	\big\la\nabla q_i(\bar x), w\big\ra =0, i\in E\cup I^+\Big\}.$$
The latter along with \eqref{tiltmod2} gives us that
$$\big\la w, \nabla_x^2\mathcal{L}(\ox, \bar\lambda)w\big\ra  \ge\frac{1}{\kk}\|w\|^2\quad \text{whenever }  	\la\nabla q_i(\ox), w\ra =0\ \mbox{as}\  i\in E\cup I^+.$$
On the other hand, by Theorem~\ref{prop53}, we have
 $$\la w, \nabla_x^2\mathcal{L}(\ox, \lambda^1)w\ra=\la w, \nabla_x^2\mathcal{L}(\ox, \lambda^2)w\ra$$ whenever $\lambda^1,\lambda^2\in \Lambda(\bar x)$
  and $\la\nabla q_i(\ox), w\ra =0\ \mbox{as}\  i\in E\cup I^+$.
  Therefore,
$$\big\la w, \nabla_x^2\mathcal{L}(\ox, \lambda)w\big\ra  \ge\frac{1}{\kk}\|w\|^2\quad \text{whenever } \ \lambda\in \Lambda(\bar x)\ \mbox{and}\  	\la\nabla q_i(\ox), w\ra =0\ \mbox{as}\  i\in E\cup I^+.$$
So, by $(i)$, the point $\bar x$ is  a  tilt-stable local minimizer of \eqref{tiltMP1} with modulus $\kappa'$, which contradicts the choice of $\kappa'.$
This shows that
\begin{equation}\label{tiltmod3}{\rm tilt}\, (\varphi, q,\ox)\leq \sup\limits_{w}\Big\{\frac{\|w\|^2} {\la w, \nabla_x^2\mathcal{L}(\ox, \bar\lambda)w\ra}\Big|\;   	\big\la\nabla q_i(\bar x), w\big\ra =0, i\in E\cup I^+\Big\}.
  \end{equation}
  By \eqref{tiltmod1} and  \eqref{tiltmod3},  we get \eqref{tiltEB}.
 The proof is complete. $\hfill\Box$

\begin{remark}  Theorem \ref{sochots} achieves the same results as  \cite[Theorem 7.7]{GM15} but under weaker, more easily verifiable assumptions. We have relaxed the CRCQ to the RCRCQ and removed the linear independence condition on the gradients of equality constraints, both of  which were requirements in \cite[Theorem 7.7]{GM15}.  The key difference in proving Theorem \ref{sochots} compared to \cite[Theorem~7.7]{GM15} lies in our use of a primal-dual approach, specifically the subgradient graphical derivative characterization of tilt-stability. In contrast, the proof of \cite[Theorem 7.7]{GM15} employed a dual-dual approach, relying on the combined second-order subdifferential characterization of tilt-stability. The latter  seemingly hinges on the linear independence condition.
 \end{remark}

Recall from  \cite[Definition 7.1]{GM15} that a  twice differentiable mapping $g: \R^m\to \R^s$ is said to be  {\it $2$-regular}  at  a given point $\ox\in\R^m$ in direction $v\in\R^m$ if for any $p\in\R^s$ there exists  $(u,w)\in\R^m\times\R^m$ such that
\begin{equation*}
\nabla g(\ox)u+[\nabla^2g(\ox)v,w]=p, \ \ \ \nabla g(\ox)w =0,
\end{equation*}
 where $[\nabla^2g(\ox)v,w]$ denotes the $s$-vector column with the entrices $\langle \nabla^2g_i(\bar x)v, w\rangle,$ $ i=1,...,s.$

For each $\ox\in\Gamma$ and  $ v\in T_\Gamma^{lin}(\ox),$   let  \begin{equation*}
I(\ox;v) := \Big\{i\in I(\ox) \ \big| \ \la\nabla q_i(\ox),v\ra=0\Big\},
\end{equation*}
  \begin{equation*} \label{1111}
\Xi(\ox;v) := \left\{z\in\R^n \ \big| \ \la\nabla q_i(\ox),z\ra+\la v,\nabla^2 q_i(\ox)v\ra\begin{cases}=0 \quad &\mbox{if}\ i\in E,\\
 \leq 0 \ &{\rm if} \ i\in I(\ox)\end{cases}\right\},
\end{equation*}
and
\begin{equation*}
 \mathcal{C}(\ox;v) := \Big\{\mathcal{C}\ \big| \ \mathcal{C} = \big\{i\in I(\ox;v)\ |\  \la\nabla q_i(\ox),z\ra+\la v,\nabla^2 q_i(\ox)v\ra = 0 \big\}\ \mbox{for some}\  z\in \Xi(\ox;v)  \Big\}.
\end{equation*}

Given a point  $\ox\in \Gamma$  and a vector $v\in K\big(\bar x, -\nabla \varphi(\bar x)\big),$ recall from \cite[Definition 7.3]{GM15} that   $\ox$  is said to be {\it nondegenerate in  direction $v$}  if the set  $\Lambda\big(\ox,-\nabla \varphi(\ox);v\big)$  is a singleton. Otherwise, one says that $\ox$ is  {\it degenerate in  direction $v$}.

A point $x$ in a convex subset $C$ of $\R^n$   is called  an {\it extreme point} of $C$ if it cannot be expressed as  $(1-t)y+tz$ for distinct $y,z\in C$  and $0<t<1$; in other words, if $x=(1-t)y+tz$ with $y,z\in C$ and $0<t<1$, then it must be that $y=z=x$ (see \cite[Section 18]{rock70}).

Theorem 7.6 in \cite{GM15} offers another pointbased characterization of tilt stability in nonlinear programming. This theorem requires, among other conditions, that for every $ v\in K\big(\ox,-\nabla \varphi(\ox)\big)\backslash\{0\}$,  one of the following must hold:
 \begin{itemize}
\item[${\bf (A)}$] either $\ox$  is   nondegenerate  in  direction $v,$
\item[${\bf (B)}$]  or for each  $\lambda\in\Lambda_{\mathcal{E}}\big(\bar x, -\nabla \varphi(\bar x); v\big)$  there exists   a maximal element $\widehat{\mathcal{C}}\in\mathcal{C}(\ox,v)$ with $I^+(\lambda)\subset\widehat{\mathcal{C}}$  such that the mapping $(q_i)_{i\in\widehat{\mathcal{C}}}$ is 2-regular at $\ox$ in direction $v$. Here,
    $$\Lambda_{\mathcal{E}}\big(\bar x, -\nabla \varphi(\bar x); v\big)=\Lambda\big(\bar x, -\nabla \varphi(\bar x); v\big)\cap \mathcal{E}\big(\bar x, -\nabla\varphi(\bar x)\big),$$
    with $\mathcal{E}\big(\bar x, -\nabla\varphi(\bar x)\big)$ being the set of all extremal points of $\Lambda\big(\bar x, -\nabla\varphi(\bar x)\big)$.
\end{itemize}

The following examples illustrate situations where Theorem \ref{sochots} is applicable, unlike both Theorem 7.6 and Theorem 7.7 from \cite{GM15}.
\begin{example}\label{ex1}{\rm Consider the following NLP in $\mathbb{R}^4$:
		\[
		 (NLP1)\quad \begin{array}{rl}  & \mbox{minimize}\quad \ \, \varphi(x):=\frac{1}{4}x_1-x_1x_2+x_3+x_3^2+\frac{1}{2}x_4^2, \quad x = (x_1, x_2, x_3, x_4)\in \R^4,\\
			& \mbox{subject to} \quad \ \,   \quad \scalebox{0.97}{$q_1(x):=x_1-x_3\leq 0, \quad q_2(x):=-x_1-x_3\leq 0, \quad q_3(x):=x_2-2x_3\leq 0, \quad$} \\
			& \hspace{1.5cm} \quad \ \, \quad q_4(x):=-x_2-2x_3\leq 0,  \quad q_5(x):=-x_1+x_2^2\leq 0, \quad q_6(x):=-x_1\leq 0, \\
			& \hspace{2.4cm} q_7(x):=-x_1+x_2^2+x_2=0.
		\end{array}
		\]		
		Put $I:=\{1, \ldots, 6\}$, $E:=\{7\}$, $q(x):=\big(q_1(x),q_2(x),...,q_7(x)\big ),$   $\Theta:= \mathbb{R}^6_{-}\times\{0\}$ and
		$$\begin{array}{rl}\Gamma&:=\{x\in\mathbb{R}^4\;|\;q_i(x)\leq 0, i\in I, \; q_i(x)= 0, i\in E\}\\
&=\{x=(x_1, x_2, x_3, x_4)\in\mathbb{R}^4\;|\; x_1 = x_2^2+x_2,\; 0\leq x_2\leq 2x_3, \; x_2^2\leq x_1\leq x_3\}.
\end{array}$$
We see that   $\ox:=(0, 0, 0, 0)\in \Gamma$,   ${I}(\ox)=I$ and
		\[\begin{array}{rl}
			& \nabla q_1(x)=(1, 0,-1, 0),  \nabla q_2(x)=(-1, 0,-1, 0),  \nabla q_3(x)=(0, 1,-2, 0),  \\
			&\nabla q_4(x)=(0, -1,-2, 0),  \nabla q_5(x)=(-1, 2x_2,0, 0),  \nabla q_6(x)=(-1, 0,0, 0),\\
			& \nabla q_7(x)=(-1, 2x_2+1,0, 0),  \nabla\varphi (x)=(\frac{1}{4}-x_2, -x_1, 1+2x_3, x_4).
		\end{array}\]
Since  $\text{rank}\{\nabla q_5(\ox), \nabla q_6(\ox)\}=1$ and  $\text{rank}\{\nabla q_5(x), \nabla q_6(x)\}=2$  for all $x=(x_1, x_2, x_3, x_4)\in\mathbb{R}^4$ with $x_2\neq 0$, the CRCQ does not hold at $\ox$.  As a result, Theorem 7.7 from \cite{GM15} cannot be used for the nonlinear program $(NLP1)$.

Next, we will verify the RCRCQ at
$\bar x$.  To do this, we choose any $J\subset I(\bar x)$ and define
$$U:=\left\{x=(x_1, x_2, x_3, x_4)\in\mathbb{R}^4|\; |x_2| <\frac{1}{4}\right\}.$$  Clearly, $U$ is a neighborhood of $\ox$. If $J=\emptyset$ then
 $$\text{rank}\{\nabla q_i(x)\, |\, i\in E\cup J\}=\text{rank}\{\nabla q_7(x)\}=1\quad \mbox{for all}\  x\in U.$$
    If $J=\{j\}$ for some $j\in I(\bar x)$ then
   $$\text{rank}\{\nabla q_i(x)\, |\, i\in E\cup J\}=\text{rank}\{\nabla q_7(x), \nabla q_j(x)\}=2\quad \mbox{for all}  \ x\in U.$$
If $J=\{5,6\}$ then
 $$\text{rank}\{\nabla q_i(x) \,| \,i\in E\cup J\}=\text{rank}\{\nabla q_5(x), \nabla q_6(x), \nabla q_7(x) \}=2\quad \mbox{for all}  \ x\in U.$$
If $J\not=\{5,6\}$ and $|J|\geq 2$ then, by a direct verification, we see that
$$\text{rank}\{\nabla q_i(x)\, |\, i\in E\cup J\}=3\quad \mbox{for all}  \ x\in U.$$
 Thus, the RCRCQ is satisfied at $\ox$.

Our aim now is to show that Theorem 7.6 from \cite{GM15} cannot be applied to the nonlinear program $(NLP1)$.
		 The set of multipliers associated with $\big(\ox, -\nabla\varphi(\ox)\big)$, denoted by $\Lambda\big(\ox, -\nabla\varphi(\ox)\big)$,  is directly calculated as:
			\begin{eqnarray*}
				\Lambda\big(\ox, -\nabla\varphi(\ox)\big)				&=&\Big\{(\lambda_1,\ldots,\lambda_7)\in\mathbb{R}^7\;|\; \lambda_i\geq 0 \text{ for } i\in I,\, \lambda_7= \lambda_4-\lambda_3,\\
				& & \hspace{2 cm}\lambda_1+\lambda_2+2\lambda_3+2\lambda_4=1,\; \lambda_1-\lambda_2+\lambda_3-\lambda_4-\lambda_5-\lambda_6=-\frac{1}{4} \Big\},
			\end{eqnarray*}
which  has more than one element.  Furthermore,
$$I^+:=\bigcup\limits_{\lambda\in \Lambda\big(\ox, -\nabla\varphi(\ox)\big)}I^+(\lambda)=I.$$
 The critical cone to $\Gamma$ at $(\ox, -\nabla\varphi(\ox))$ is computed as follows:
			\begin{eqnarray*}
				K(\ox, -\nabla\varphi(\ox))&=&T_{\Gamma}(\ox)\cap \{\nabla\varphi(\ox)\}^{\bot}\\
				&=& \{u\in\mathbb{R}^4\;|\;  \la\nabla q_i(\bar x), u\ra\leqslant 0, i\in I, \la\nabla q_7(\bar x), u\ra= 0 \}\cap \{\nabla\varphi(\ox)\}^{\bot}\\
				&=&\{(0,0, 0)\}\times \mathbb{R}.
			\end{eqnarray*}
	 Take any  $v\in K\big(\ox, -\nabla\varphi(\ox)\big).$  Then $v$ satisfies \eqref{msid} with $x=\bar x$ and $x^*=-\nabla\varphi(\bar x)$.  On the other hand,  RCRCQ holds at $\ox$. Therefore, by Theorem \ref{prop53},
$$\Lambda\big(\ox, -\nabla\varphi(\ox); v\big) = \Lambda\big(\ox, -\nabla\varphi(\ox)\big),$$
which is not a singleton. This justifies that $\ox$ degenerates in the critical directions $v$.
Let $\mathcal{E}\big(\ox, -\nabla\varphi(\ox)\big)$ be the set of all extreme points of the  multiplier set $\Lambda\big(\ox, -\nabla\varphi(\ox)\big)$, and let $\lambda\in \Lambda\big(\ox, -\nabla\varphi(\ox)\big)$. Then $\lambda\in\mathcal{E}\big(\ox, -\nabla\varphi(\ox)\big)$ if and only if
 $\{\nabla q_i(\bar x)| i\in E\cup I^+(\lambda)\}$ is linearly independent. Suppose that $\lambda\in\mathcal{E}\big(\ox, -\nabla\varphi(\ox)\big).$ Then $\{\nabla q_i(\bar x)| i\in E\cup I^+(\lambda)\}$ is linearly independent. This along with  $I^+(\lambda)\subset I$ and  $\mbox{rank}\left\lbrace \nabla q_i(\bar x) \right\rbrace_{i\in E\cup I}\leq 3$ tells us that  $|I^+(\lambda)|\leq 2.$ Thus,  there are at least four indices $i\in I$ such that  $\lambda_i=0$.
\par
  {\it Case 1:}  $\lambda_i=0$ for   $i\in\{3, 4, 5, 6\}$. Then $\lambda_1=3/8>0$, $\lambda_2=5/8>0$, and $\lambda_7=0$.
\par
{\it Case 2:} $\lambda_i=0$ for   $i\in\{2, 4, 5, 6\}$. Then $\lambda_1=-2$, $\lambda_3=5/4>0$, and $\lambda_7=-5/4$. This is impossible since $\lambda_1<0.$
\par
{\it Case 3:}  $\lambda_i=0$ for  $i\in\{2, 3, 5, 6\}$. Then $\lambda_1=1/6>0$,  $\lambda_4=5/12>0$,  and $\lambda_7=5/12$.
\par
{\it Case 4:}  $\lambda_i=0$ for  $ i\in\{2, 3, 4, 6\}$.  Then $\lambda_1=1>0$,  $\lambda_5=5/4>0$,  and $\lambda_7=0$.
\par{\it Case 5:}  $\lambda_i=0$ for  $ i\in\{2, 3, 4, 5\}$.  Then $\lambda_1=1>0$,  $\lambda_6=5/4>0$,  and $\lambda_7=0$.
\par{\it Case 6:}  $\lambda_i=0$ for  $ i\in\{1,  4, 5, 6\}$.  Then $\lambda_2=1/2>0$,  $\lambda_3=1/4>0$,  and $\lambda_7=-1/4$.
\par{\it Case 7:}  $\lambda_i=0$ for  $ i\in\{1,  3, 5, 6\}$.  Then $\lambda_2=-1/2$,  $\lambda_4=3/4>0$,  and $\lambda_7=3/4$. This is impossible since $\lambda_2<0.$
\par{\it Case 8:}  $\lambda_i=0$ for  $ i\in\{1,  3, 4, 6\}$.  Then $\lambda_2=1$,  $\lambda_5=-3/4$,  and $\lambda_7=0$. This is impossible since $\lambda_5<0.$
\par{\it Case 9:}  $\lambda_i=0$ for  $ i\in\{1,  3, 4, 5\}$.  Then $\lambda_2=1$,  $\lambda_6=-3/4$,  and $\lambda_7=0$. This is impossible since $\lambda_6<0.$
\par{\it Case 10:}  $\lambda_i=0$ for  $ i\in\{1, 2, 5, 6\}$.  Then $\lambda_3=1/8>0$,  $\lambda_4=3/8>0$,  and $\lambda_7=1/4$.
\par{\it Case 11:}  $\lambda_i=0$ for  $ i\in\{1, 2, 4, 6\}$.  Then $\lambda_3=1/2>0$,  $\lambda_5=3/4>0$,  and $\lambda_7=-1/2$.
\par {\it Case 12:}  $\lambda_i=0$ for  $ i\in\{1, 2, 4, 5\}$. Then $\lambda_3=1/2>0$,  $\lambda_6=3/4>0$,  and $\lambda_7=-1/2$.
\par {\it Case 13:}  $\lambda_i=0$ for  $ i\in\{1, 2, 3, 6\}$. Then $\lambda_4=1/2>0$,  $\lambda_5=-1/4$,  and $\lambda_7=1/2$.
This is impossible since $\lambda_5<0.$
\par {\it Case 14:}  $\lambda_i=0$ for  $ i\in\{1, 2, 3, 5\}$. Then $\lambda_4=1/2>0$,  $\lambda_6=-1/4$,  and $\lambda_7=1/2$.
This is impossible since $\lambda_6<0.$
\par {\it Case 15:}  $\lambda_i=0$ for  $ i\in\{1, 2, 3, 4\}$.
This is impossible since $\lambda_1+\lambda_2+2\lambda_3+2\lambda_4=1.$
\par\noindent
Therefore,
		\begin{equation}\label{etrEq}
\begin{array}{rl}
				\mathcal{E}\big(\ox, -\nabla\varphi(\ox)\big)\subset& \bigg\{ \left(\frac{3}{8}, \frac{5}{8}, 0, 0, 0, 0, 0 \right), \; \left(\frac{1}{6},  0, 0, \frac{5}{12}, 0, 0, \frac{5}{12} \right), \; \left(1,  0, 0, 0, \frac{5}{4}, 0, 0 \right), \\
&\  \left(1,  0, 0, 0,  0, \frac{5}{4}, 0 \right), \; \left(0, \frac{1}{2}, \frac{1}{4}, 0, 0,  0, -\frac{1}{4} \right), \; \left(0, 0, \frac{1}{8}, \frac{3}{8}, 0, 0, \frac{1}{4} \right),
\\
&			\left(0, 0, \frac{1}{2}, 0,  \frac{3}{4}, 0, -\frac{1}{2} \right), \;	\left(0, 0, \frac{1}{2}, 0, 0,  \frac{3}{4}, -\frac{1}{2} \right)		
			\bigg\}.
\end{array}
		\end{equation}
Conversely, let $\lambda$  be any element from the right-hand side set of \eqref{etrEq}. By a direct verification, we see that $\{\nabla q_i(\bar x)| i\in E\cup I^+(\lambda)\}$	 is linearly independent. This shows that $\lambda\in \mathcal{E}(\ox, -\nabla\varphi(\ox)\big)$. Thus,  by \eqref{etrEq}, we have
	\begin{equation*}
\begin{array}{rl}
				\mathcal{E}\big(\ox, -\nabla\varphi(\ox)\big)=& \bigg\{ \left(\frac{3}{8}, \frac{5}{8}, 0, 0, 0, 0, 0 \right), \; \left(\frac{1}{6},  0, 0, \frac{5}{12}, 0, 0, \frac{5}{12} \right), \; \left(1,  0, 0, 0, \frac{5}{4}, 0, 0 \right), \\
&\  \left(1,  0, 0, 0,  0, \frac{5}{4}, 0 \right), \; \left(0, \frac{1}{2}, \frac{1}{4}, 0, 0,  0, -\frac{1}{4} \right), \; \left(0, 0, \frac{1}{8}, \frac{3}{8}, 0, 0, \frac{1}{4} \right),
\\
&			\left(0, 0, \frac{1}{2}, 0,  \frac{3}{4}, 0, -\frac{1}{2} \right), \;	\left(0, 0, \frac{1}{2}, 0, 0,  \frac{3}{4}, -\frac{1}{2} \right)		
			\bigg\}.
\end{array}
		\end{equation*}  	
 Furthermore, the extreme multiplier set $\Lambda_{\mathcal{E}}\big(\ox, -\nabla\varphi(\ox); v\big)$ in  direction $v$ is computed as follows:
		\begin{eqnarray*}
			\Lambda_{\mathcal{E}}\big(\ox, -\nabla\varphi(\ox); v\big)&=& \Lambda\big(\ox, -\nabla\varphi(\ox); v\big) \cap \mathcal{E}\big(\ox, -\nabla\varphi(\ox)\big)\\
			&=&\Lambda\big(\ox, -\nabla\varphi(\ox)\big) \cap \mathcal{E}\big(\ox, -\nabla\varphi(\ox)\big)\\
			&=&\mathcal{E}\big(\ox, -\nabla\varphi(\ox)\big);
		\end{eqnarray*}
see \cite[p. 2100]{GM15}.
Note  that $v_1=v_2=v_3=0,$ and  $${I}(\ox; v):=\left\{i\in{I}(\ox)|\; \la\nabla q_i(\ox),v\ra =0\right\} = {I}(\ox) = I,$$ and   $\nabla^2 q_i(x)=0$ for  $i= 1, 2, 3, 4, 6$, and $\nabla^2q_5(x)=\nabla^2q_7(x)= \text{diag}(0, 2, 0, 0)$, where $\text{diag}(a_1, a_2,$ $\ldots, a_n)$ denotes the square diagonal matrix with the main diagonal entries $a_1, a_2,$ $\ldots, a_n$.  We get that
		  \begin{eqnarray*}
		  	 \Xi(\ox; v)&:=&\left\lbrace z\in \mathbb{R}^4\bigg|\; \la\nabla q_i(\ox),z\ra + \la v, \nabla^2q_i(\ox)v\ra  \begin{cases}
		  		=0 &\text{ for } i\in E\\
		  		\leq 0 &\text{ for } i\in {I}(\ox)
		  	\end{cases} \right\rbrace \\
	  	&=& \Big\{ (z_1, z_2, z_3, z_4)\in \mathbb{R}^4\;|\; z_1-z_3\leqslant 0, -z_1-z_3\leqslant 0, z_2-2z_3\leqslant 0,  \\
	  	& &\hspace{3cm}-z_2-2z_3\leqslant 0, -z_1+2v_2^2\leqslant 0, -z_1\leqslant 0, -z_1+z_2+2v_2^2=0 \Big\} \\
	  	&=& T_{\Gamma}(\ox),
		  \end{eqnarray*}
and
\begin{eqnarray*}
	\mathcal{C}(\ox; v)&:=\left\lbrace \mathcal{C} \subset{I}(\ox; v)\Big|\, \exists z\in \Xi(\ox; v)\text{ with }
	\mathcal{C}=\left\lbrace i\in {I}(\ox)|\; \la\nabla q_i(\ox),z\ra + \la v, \nabla^2q_i(\ox)v\ra =0 \right\rbrace \right\rbrace\\
	& =\left\lbrace \mathcal{C} \subset {I}(\ox)\Big| \, \exists z\in T_{\Gamma}(\ox)\text{ with }
	\mathcal{C}=\left\lbrace i\in {I}(\ox)|\; \la\nabla q_i(\ox),z\ra =0 \right\rbrace \right\rbrace.
\end{eqnarray*}
Since  $\la\nabla q_i(\ox),z\ra =0$  for all  $i\in I(\bar x)$  and $z=0_{\mathbb{R}^4}\in \Xi(\ox; v)$, it follows that
$I(\bar x)\in \mathcal{C}(\ox; v).$ This implies that  $I(\bar x)$ is the unique maximal element in $\mathcal{C}(\ox; v),$ and  for every $\lambda\in\Lambda_{\mathcal{E}}\big(\ox, -\nabla\varphi(\ox); v\big)$  the set $\widehat{\mathcal{C}}:=I(\bar x)$  is 	a maximal element in $\mathcal{C}(\ox; v)$ satisfying $I^+(\lambda) \subset \widehat{\mathcal{C}}$.
We assert that the narrow active constraint mapping $(q_i)_{i\in E\cup \widehat{\mathcal{C}}}$ is not 2-regular at $\ox$ in the critical direction $v\neq 0$.
 Let us assume, for contradiction, that our claim is false. Then, the mapping  $q: \mathbb{R}^4\to \mathbb{R}^7$, defined as $q(x)=\big(q_1(x), q_2(x),...,q_7(x)\big)$,   is  2-regular at $\ox$ in the critical direction $v\neq 0$. By the definition of 2-regularity, for every $p\in \R^7$ the system
	\begin{eqnarray}\label{sys1}
		\nabla q(\ox)u + [\nabla^2q(\ox)v, w] = p,\; \nabla q(\bar x) w=0
	\end{eqnarray}
	admits a solution $(u, w)\in \mathbb{R}^4\times \mathbb{R}^4$.  Here, the symbol $[\nabla^2q(\ox)v, w]$
represents the	7-vector column with  entries $\la \nabla^2q_i(\ox)v, w\ra$ for  $i = 1,\ldots,7$.
 On the other hand, since $ [\nabla^2q(\ox)v, w]=0$, we can equivalently express \eqref{sys1} as:
		$$\left\{
	\begin{array}{l}
		\nabla q(\ox)u  = p\\
		\nabla q(\bar x) w=0.			
	\end{array}
	\right.$$
Therefore, the linear mapping $\nabla q(\ox): \R^4\to \R^7$ is surjective, which is a contradiction because ${\rm dim}\nabla q(\ox)(\R^4)\leq 4<7$.
This proves that the narrow active constraint mapping $(q_i)_{i\in E\cup \widehat{\mathcal{C}}}$ is not 2-regular at $\ox$ in the critical direction $v$.
Consequently, since the set $K\big(\ox, -\nabla\varphi(\ox)\big)\backslash\{0\}$ is nonempty,  and for every $v$ in this set, neither condition ({\bf A}) nor condition ({\bf B}) is satisfied, we can conclude that Theorem 7.6 from \cite{GM15} cannot be applied to the nonlinear program $(NLP1)$.

Our final step is to apply Theorem \ref{sochots} to demonstrate that
$\bar x$   is a tilt-stable local minimizer of $(NLP1)$.  The Lagrange function is given by
			\begin{eqnarray*}
				\mathcal{L}(x, \lambda)&=&\varphi(x)+\lambda_1q_1(x)+\ldots+\lambda_7q_7(x)\quad \text{for all } (x, \lambda)\in \mathbb{R}^4\times\mathbb{R}^7.
			\end{eqnarray*}
			Thus, we have
			\[\nabla_x\mathcal{L}(x, \lambda) =\begin{bmatrix}
				\frac{1}{4}-x_2+\lambda_1-\lambda_2-\lambda_5-\lambda_6-\lambda_7\\
				-x_1+\lambda_3-\lambda_4+2\lambda_5x_2+2\lambda_7x_2+\lambda_7\\
				1+2x_3-\lambda_1-\lambda_2-2\lambda_3-2\lambda_4\\
				x_4
			\end{bmatrix}\  \mbox{and}\,\   \nabla^2_x\mathcal{L}(x, \lambda)=\begin{bmatrix}
				0&-1&0&\;0\\ -1&\;\;2\lambda_5+2\lambda_7\;\;&0&\; 0  \\0&0&2&\; 0\\0&0&0&\;1				
			\end{bmatrix}.\]
Since  $\Lambda\big(\ox, -\nabla\varphi(\ox)\big)\not=\emptyset$, we know that $\bar x$ is a stationary point of $(NPL1)$. Furthermore,  we have $$	\big\{w\in \mathbb{R}^4\;|\; \la\nabla q_i(\ox), w\ra = 0, i\in E\cup I^+\big\} = \{(0,0, 0)\}\times \mathbb{R}.$$ Thus, for every  $\lambda\in	\Lambda\big(\ox, -\nabla\varphi(\ox)\big)$ and for every nonzero $w=(w_1, w_2, w_3, w_4)\in \mathbb{R}^4$ with $\la\nabla q_i(\ox), w\ra = 0$ as $i\in E\cup I^+$
we see that  $w_1=w_2=w_3=0,$  $w_4\neq 0$ , and
		\[\la w, \nabla_x^2\mathcal{L}(\ox, \lambda)w\ra=w_4^2 >0.\]
		Therefore,  by Theorem \ref{sochots},  $\bar{x}$ is a tilt-stable minimizer of $(NLP1),$ and   the exact bound of tilt stability is  calculated by
$$
	{\rm tilt}\, (\varphi, q,\ox)=\sup\limits_{w}\left\{\frac{\|w\|^2} {\la w, \nabla_x^2\mathcal{L}(\ox, \bar\lambda)w\ra}\Big|\;   	\big\la\nabla q_i(\bar x), w\big\ra =0, i\in E\cup I^+\right\}=\sup\limits_{w_4\in \R}\frac{w_4^2} {w_4^2}=1,
$$
where $\bar\lambda$ is taken arbitrarily from the multiplier set $\Lambda(\bar x):=\Lambda\big(\ox, -\nabla\varphi(\ox)\big)$.

	}
\end{example}

\begin{example}\label{ex2}{\rm Consider the following NLP in $\mathbb{R}^4$:		
\[
 (NLP2) \quad \begin{array}{rl}  & \mbox{minimize}\quad \ \, \varphi(x):=\frac{1}{4}x_1-x_1x_2+x_3+x_3^2+\frac{1}{12}x_4^4, \quad x = (x_1, x_2, x_3, x_4)\in\mathbb{R}^4,\\
	& \mbox{subject to}  \ \,   \quad q_1(x):=x_1-x_3\leq 0, \, q_2(x):=-x_1-x_3\leq 0, \, q_3(x):=x_2-2x_3\leq 0, \quad \\
	& \hspace{1.5cm}  \ \, \quad q_4(x):=-x_2-2x_3\leq 0,  \, q_5(x):=-x_1+x_2^2\leq 0, \, q_6(x):=-x_1\leq 0, \\
	& \hspace{2.05cm} q_7(x):=-x_1+x_2^2+x_2=0.
\end{array}
\]	
This nonlinear program  is a variation of the previous one, with the only change being the objective function where the term $\frac{1}{2}x_4^2$ is replaced by $\frac{1}{12}x_4^4$.
 Thus, the first-order derivatives at $\bar x=0_{\R^4}$ remain unchanged, and as a result, the assumptions of Theorem \ref{sochots} are still met at this point.  However, the characterizations of tilt stability for NLPs under CRCQ \cite[Theorem 7.7]{GM15} and under either nondegeneracy or 2-regularity \cite[Theorem 7.6]{GM15} are not applicable here.
The Hessian of the Lagrangian is given by
			\[  \nabla^2_x\mathcal{L}(x, \lambda)=\begin{bmatrix}
			0&-1&0&\;0\\ -1&\;\;2\lambda_5+2\lambda_7\;\;&0&\; 0  \\0&0&2&\; 0\\0&0&0&\;x_4^2				
		\end{bmatrix}.\]
Thus, choosing  $w=(0, 0, 0, w_4)\in\mathbb{R}^4$ with $w_4\neq 0$, we find that $\la\nabla q_i(\ox), w\ra = 0$ as $i\in E\cup I^+$,  and
			$\la w, \nabla_x^2\mathcal{L}(\ox, \lambda)w\ra=0$
	 for every $\lambda\in	\Lambda\big(\ox, -\nabla\varphi(\ox)\big)$.	 By Theorem \ref{sochots}, $\ox$ is not a tilt-stable minimizer of $(NLP2).$
	}
\end{example}
 \section{Concluding Remarks}

In this paper, we have utilized the neighborhood primal-dual approach to tilt stability and extended calculus rules for the subgradient graphical derivative to characterize tilt-stable local minimizers and calculate  the exact bound of tilt stability for the NLPs under RCRCQ.
Our work extends \cite[Theorem~7.7]{GM15}  by relaxing the constraint qualification and removing the linear independence condition for equality constraint gradients. This significantly contributes to our understanding of tilt stability for NLPs.
Based on our results, particularly the simplified graphical derivative formula for normal cone mappings (Theorem 3.1) and the pointbased second-order characterizations of tilt-stable minimizers for the NLPs under RCRCQ (Theorem 4.1), two key questions arise for future research:

1. Can we precisely compute the limiting coderivative of normal cone mappings associated with the feasible sets of NLPs under RCRCQ?

2. Can the neighborhood primal-dual approach be applied to the study of full stability, in the sense of Levy, Poliquin, and Rockafellar \cite{LPR00}, for NLPs under RCRCQ?

\vspace{0.3cm}
{\bf Acknowledgments.} This research is funded by Funds for Science and Technology Development of the University of Danang under project number B2023-DN03-08.

\end{document}